\documentclass{article}

\usepackage[utf8]{inputenc}
\usepackage[T1]{fontenc}
\usepackage[english]{babel}
\usepackage{lmodern}
\usepackage{calrsfs}
\usepackage{color}
\usepackage{amsmath}
\usepackage{amsthm}
\usepackage{mathtools, bm}
\usepackage{amssymb, bm}
\usepackage{array}
\usepackage{enumerate}
\usepackage{sectsty}
\usepackage{graphicx}
\usepackage{float}
\usepackage{cancel}
\usepackage{url}
\usepackage{nomencl}
\makenomenclature
\usepackage{makeidx}
\makeindex

\usepackage{here}
\usepackage{geometry}
\usepackage[pdftex]{hyperref}
\usepackage{pifont}
\usepackage{dsfont}

\usepackage{pgf}
\usepackage{tikz}
\usetikzlibrary{matrix,arrows,automata}
\usetikzlibrary{shapes}
\usetikzlibrary{calc}
\usetikzlibrary{arrows}
\usetikzlibrary{decorations.pathreplacing,decorations.markings}

\chapterfont{\centering}
\sectionfont{\centering}
\subsectionfont{}

\graphicspath{{C:/Users/gwlad/Desktop/Gwladys/These/article/Equations_Schroder/}}

\newtheoremstyle{monstyledem} 
{8pt}                    
{8pt}                    
{\normalfont}                   
{}                           
{\bf}                   
{\newline}                          
{.5em}                       
{}  

\newtheoremstyle{monstyle} 
{8pt}                    
{8pt}                    
{\itshape}                   
{}                           
{\bf}                   
{\newline}                          
{.5em}                       
{}  

\theoremstyle{monstyle}

\newtheorem{thm_en}{Theorem}[section]

\theoremstyle{monstyledem}

\newcommand\blankfootnote[1]{%
	\let\thefootnote\relax\footnotetext{#1}%
	\let\thefootnote\svthefootnote%
}
\newcommand \N{\mathbb{N}}
\newcommand \R{\mathbb{R}}
\newcommand \C{\mathbb{C}}
\newcommand \ES{E. Schröder }
\newcommand \Beck{P.-G. Becker }
\newcommand \Berg{W. Bergweiler }
\newcommand \Ritt{J. F. Ritt }

\title{A survey on the hypertranscendence of the solutions of the Schröder's, Böttcher's and Abel's equations}
\author{Gwladys Fernandes}
\date{}

\begin{document}
	\maketitle
	
\textbf{Abstract.} In 1994, P.-G. Becker and W. Bergweiler \cite{Becker_Bergweiler} listed all the differentially algebraic solutions of three famous functional equations: the Schröder's, Böttcher's and Abel's equations. The proof of this theorem combines various domains of mathematics. This goes from the theory of iteration, which gave birth to these equations, to the algebro-differential notion of coherent families developed by M. Boshernitzan and L. A. Rubel. This survey is an excursion into the history of these equations, in order to enlighten the different pieces of mathematics they bring together and how these parts fit into the result of P.-G. Becker and W. Bergweiler.

\section{Introduction}

This survey is about three famous functional equations, named after the mathematicians Ernst Schröder, Lucien Böttcher and Niels Henrik Abel. These equations are linked to the local study of the iteration of a rational function $R(z)\in\C(z)$ with complex coefficients around a \textbf{fixed point} $\alpha$. That is a point of $\C\cup\{\infty\}$ such that $R(\alpha)=\alpha$. We can assume, without any loss of generality, that $\alpha=0$ (see Section \ref{sec_birth} and \eqref{sol_eq_conj}). Moreover, to avoid the trivial case of Möbius transformations (see \eqref{def_Mobius} for the definition), which is well-understood, we assume that the \textbf{degree} of $R(z)$, that is the maximum of the degrees of the coprime polynomials on its numerator and denominator, is at least $2$. Then, the equations we are interested in are the following:

\begin{enumerate}

\item
Let $s=R'(0)$. If $s\neq 0$, then, the Schröder's equation is:
\begin{equation*}
\label{S_eq}
f(s z)=R(f(z)), \tag{S}
\end{equation*}	

\item

If $R(z)=\sum_{n=d}^{+\infty}a_{n}z^{n}$, where $d\geq 2$, then the Böttcher's equation is:

\begin{equation*}
\label{B_eq}
f(z^{d})=R(f(z)), \tag{B}
\end{equation*}
\item

The Abel's equation is:

\begin{equation*}
\label{A_eq}
f(R(z))=f(z)+1. \tag{A}
\end{equation*}

\end{enumerate}

In 1994, P.-G. Becker and W. Bergweiler \cite{Becker_Bergweiler} listed all the differentially algebraic solutions of Equations $(S)$, $(B)$ and $(A)$. The aim of this survey is to present the proof of this result as a testimony of the richness of the interactions these three equations centralize between various areas of mathematics and how the two authors combine them to get their beautiful result. For more historical details on the theory of iteration, besides the ones given below, we refer to \cite{Audin} and \cite{Alexander}.
\medskip

Before diving into the history of these equations, let us remind some definitions. First, a formal power series $f(z)$ with coefficients in the complex plane $\C$ is said to be \textbf{differentially algebraic} over $\C(z)$ if there exists a non-zero polynomial $P(z,X_{0},\ldots,X_{n})$ with coefficients in $\C$ such that:
\begin{equation}
\label{def_DA}
P(z,f(z),f'(z),\ldots, f^{(n)}(z))=0,
\end{equation} 
where $f^{(n)}(z)$ is the $n$-th derivative of $f$ with respect to $z$. We say that $f$ is \textbf{differentially transcendental} or \textbf{hypertranscendental} over $\C(z)$ if it is not differentially algebraic over $\C(z)$. This notion of differential algebraicity generalises the one of algebraicity. Indeed, a formal power series $f$ with coefficients in $\C$ is said to be \textbf{algebraic} over $\C(z)$ if there exists a non-zero polynomial $P(z,X)$ with coefficients in $\C$ such that:
\begin{equation}
P(z,f(z))=0,
\end{equation}
and it is said to be \textbf{transcendental} over $\C(z)$ otherwise. Moreover, formal power series $f_{1}(z),\ldots,f_{n}(z)$ with coefficients in $\C$ are said to be \textbf{algebraically dependent} over $\C(z)$ if there exists a non-zero polynomial $P(z,X_{0},\ldots,X_{n})$ with coefficients in $\C$ such that:
\begin{equation}
P(z,f_{1}(z),\ldots,f_{n}(z))=0.
\end{equation} 
If they are not algebraically dependent over $\C(z)$, we say that these functions are \textbf{algebraically independent} over $\C(z)$. Thus, a hypertranscendental function is transcendental, and all its derivatives are algebraically independent over $\C(z)$.
\medskip

Coming back to our topic, the three equations $(S), (B)$ and $(A)$ are introduced for the need of the iteration theory, for which Newton's method is a famous representative, and the research of fixed points of a rational fraction. They reach their zenith with the development of the theory of P. Fatou and G. Julia around 1918, which divides the complex plane into different domains according to the local behaviour of the iterates of a rational fraction, in terms of convergence or divergence to a fixed point. Thus, the rather algebraic problem of iteration of a rational fraction and the determination of its fixed points is linked to the rather analytic theory developed by P. Fatou and G. Julia. The interface between these different domains of mathematics grows in the following years with the work of J. F. Ritt \cite{Ritt} in the 20's, who furnishes the list of all the differentially algebraic solutions of the Schröder's equation near a repelling fixed point (see the definition in Section \ref{sec_birth}). This is a first step toward the result of  P.-G. Becker and W. Bergweiler \cite{Becker_Bergweiler} we are interested in (Theorem \ref{thm_BB} of the present paper). After that, the interest for these equations wanes for almost sixty years, with a renewal of enthusiasm in the 80's. Indeed, then, the sets of Fatou and Julia share connexions with dynamical systems and fractals which are very prolific areas at that time. From the algebraic point of view, this surge of interest is visible in 1995 in the result of P.-G. Becker and W. Bergweiler \cite{Becker_Bergweiler}, which completes the previous result of J. F. Ritt on the Schröder's equation. The two authors also provide a partial result of Theorem \ref{thm_BB} one year earlier in \cite{Becker_Bergweiler_0}. Indeed, they give the list of all the algebraic solutions of the Böttcher's equation $(B)$. In this paper, we present their accomplished result \cite{Becker_Bergweiler}, reproduced in Theorem \ref{thm_BB}. This statement provides the list of all the differentially algebraic solutions of the Schröder's, Böttcher's and Abel's equations. The proof of P.-G. Becker and W. Bergweiler relies on the results in \cite{Ritt, Becker_Bergweiler_0}, a theorem on coherent families developed by M. Boshernitzan and L. A. Rubel \cite{Boshernitzan_Rubel}, and the theory of P. Fatou and G. Julia \cite{Fatou_1,Fatou_2,Fatou_3,Julia}.

\section{Schröder's, Böttcher's and Abel's equations: their deep interactions with the iteration theory through history}

\subsection{The birth of the Schröder's, Böttcher's and Abel's equations} 
\label{sec_birth}

The Schröder's equation appears in 1870 in a paper of the same author \cite{Schroder_1}, in link with the strong interest of the author in Newton's method. This consists in finding approximations of real roots of a real-valued function $f(x)$. The idea is the following. Take a real number $x_0$ and define the following sequence for $n\geq 0$: 
\begin{equation}
\label{Newton}
x_{n+1}=x_{n}-\frac{f(x_{n})}{f'(x_{n})}.
\end{equation}
Graphically speaking, the number $x_{n+1}$ is the abscissa of the intersection point between the tangent of the curve of $f(x)$ at $x=x_{n}$ and the $x$-axis. Then, the general principle is that, if one takes $x_{0}$ close enough to a root of $f(x)$, then the sequence $(x_{n})$ may converge to this root. As mentioned in \cite{Alexander}, this problem seems to be one of the oldest processes of iteration in the history of mathematics and we can find tracks of this in ancient Babylone or in the Arab world of the twelfth century.

We shall point out that I. Newton did not present his method with the expression of the sequence \eqref{Newton} but with an equivalent approach based on algebra rather than calculus. An other equivalent formulation is then made by J. Raphson, again without the use of calculus. It is finally T. Simpson who introduces the closest procedure to \eqref{Newton},  formally defined in this precise form by J. Fourier.
\medskip

The contribution of E. Schröder to the further development and generalization of Newton's method results from his idea of turning the discrete problem of the convergence of the sequence \eqref{Newton} into the iteration problem of the function: 
\begin{equation}
\label{Newton_frac}
R(x)=x-\frac{f(x)}{f'(x)}.
\end{equation} 
With this point of view, a zero $\alpha$ of $f(x)$ becomes a fixed point of $R(x)$. This also allows the author to extend Newton's method to the complex plane and to the search of complex fixed points of $R(x)$, or complex zeros of $f(x)$. We will be mainly interested in the case where $R(z)$ is a rational fraction with coefficients in $\C$ of degree at least two. Let $R^{n}(z)$ denote the $n$-th iterate of $R(z)$. We remark that if we start with a point $z_{0}$ and if $R^{n}(z_{0})$ converges to a point $\alpha$, then $R(\alpha)=\alpha$. That is, $\alpha$ is a fixed point of $R(z)$. This explains why such points are important in the theory of iteration. Now, by Taylor formula, in a neighbourhood of the fixed point $\alpha$, the value $|R(z)-\alpha|$ can be approximated by $|R'(\alpha)||z-\alpha|$. Hence, we guess that the behaviour of the sequence $R^{n}(z)$ is not the same depending whether $|R'(\alpha)|<1$ or $|R'(\alpha)|>1$ for example. Indeed, in the first case, the sequence converges to $\alpha$ (it is the fixed point theorem of \ES and G. Koenigs mentioned below), and in the second, the only way for this sequence to converge to $\alpha$ is that $R^{N}(z)=\alpha$ for some $N\in\mathbb{N}$. That is why we distinguish the following categories of fixed points $\alpha$ of $R$, according to the value of $|R'(\alpha)|$. We say that a fixed point $\alpha$ is:
\begin{enumerate}
	\item 
	\textbf{attracting} if $0<|R'(\alpha)|<1$,
	\item 
	\textbf{super-attracting} if $|R'(\alpha)|=0$,
	\item 
	\textbf{repelling} if $|R'(\alpha)|>1$,
	\item 
	\textbf{rationally indifferent} if $R'(\alpha)$ is a root of unity,
	\item 
	\textbf{irrationally indifferent} if $|R'(\alpha)|=1$ and if $R'(\alpha)$ is not a root of unity.	
\end{enumerate}

E. Schröder establishes (even if his proof is not rigorously explained) the following fixed point theorem in \cite{Schroder_1}: if $R(z)$ is an analytic function in a neighbourhood of an attracting fixed point $\alpha$ of $R$, then, there exists a neighbourhood $D$ of $\alpha$ such that:
\begin{equation*}
\lim_{n\to +\infty}R^{n}(z)=\alpha, \forall z\in D.
\end{equation*}

Note that G. Koenigs provides a complete proof of this theorem in \cite{Koenigs_1}. This result can be obtained by applying the Taylor formula to $R(z)$ near $z=\alpha$. As $\alpha$ is attracting, we find the existence of a neighbourhood $D$ of $\alpha$, and a real number $\epsilon$ such that $0<\epsilon<1$ and :
\begin{equation*}
|R(z)-\alpha|<\epsilon|z-\alpha|, \text{ } \forall z\in D.
\end{equation*}
We deduce that $R(D)\subset D$. By induction on $n$, we find that $R^{n}(D)\subset D$, for every integer $n\geq 1$ and the announced convergence in $D$.

Now, let us assume that the rational fraction $R(z)$ of the Newton's method \eqref{Newton_frac} satisfies the assumption of the fixed point theorem. Then, the aim of \ES is to generalize Newton's method by finding rational fractions $\phi(z)$ distinct from $R(z)$, such that $\phi(\alpha)=\alpha$ and such that $\phi^{n}(z)$ converges to $\alpha$ when $z$ is close enough to $\alpha$, in order to improve the rate of convergence of $R^{n}(z)$. To do so, it is important to understand the iterates of a rational fraction. But these are in general difficult to compute. That is why \ES thinks about finding these $\phi(z)$ such that their iterates are easy to compute, while keeping track of the initial rational fraction $R(z)$. He solves this problem by the use of conjugation. We say that two rational fractions $\phi$ and $R$ are \textbf{conjugated} if there exists an invertible function $F(z)$ such that:
\begin{equation}
\label{conj_R}
R=F^{-1}\phi F,
\end{equation} 
where functional composition is multiplicatively written.
Particularly interesting conjugations are those with a \textbf{Möbius transformation} $F$, that is a rational fraction of the following form:
\begin{equation}
\label{def_Mobius}
F(z)=\frac{az+b}{cz+d},
\end{equation}
where $a,b,c,d\in\C$ and $ad-bc\neq 0$. They are of degree one and stable under composition.

Let us stress that the conjugation preserves the notions of fixed points and iteration. Indeed, if $\alpha$ is a fixed point of $R(z)$, of one of the five kinds defined before (attracting, super-attracting, repelling, rationally indifferent or irrationally indifferent) then $F(\alpha)$ is a fixed point of $\phi(z)$ \textbf{of the same kind}, and for every $n\in\N$:
\begin{equation*}
R^{n}=F^{-1}\phi^{n} F.
\end{equation*}

The notion of fixed points is crucial in the study of equations $(S)$, $(B)$ and $(A)$. If a rational fraction $R$ admits a fixed point $\alpha$, we let $S_{R,\alpha}$, $B_{R,\alpha}$, $A_{R,\alpha}$ denote respectively the Schröder's, Böttcher's and Abel's equation $(S), (B)$ and $(A)$ associated to the rational fraction $R$. The mention of $\alpha$ means that we are interested in the behaviour of the iterations of $R(z)$ in a neighbourhood of the fixed point $\alpha$. 
\medskip

Furthermore, the conjugation respects the solutions of Equations $(S), (B)$ and $(A)$. Indeed, if $R(z)$ is a rational fraction which admits a fixed point $\alpha$, then we have the following:

\begin{enumerate}
\item 
\[f \text{ is a solution of } S_{R,\alpha} \Leftrightarrow gf \text{ is a solution of } S_{gRg^{-1},g(\alpha)}.\]
\item 
\begin{equation}
\label{sol_eq_conj}
f \text{ is a solution of } B_{R,\alpha} \Leftrightarrow gf \text{ is a solution of } B_{gRg^{-1},g(\alpha)}.
\end{equation}
\item 
\[f \text{ is a solution of } A_{R,\alpha} \Leftrightarrow fg^{-1} \text{ is a solution of } A_{gRg^{-1},g(\alpha)}.\]
\end{enumerate}

Let us go back to the approach of E. Schröder to conjugate the rational fraction of Newton's method to find an easier form. The choice of $\phi(z)=sz$ in \eqref{conj_R}, with a certain $s\in \C^{*}$ leads to the Schröder's equation:
\begin{equation}
\tag{$S_{0}$}
F(R(z))=sF(z),
\end{equation}
in which the unknown is the function $F(z)$.

Note that this is slightly different from Equation $(S)$, with the following link: if $f$ is an invertible solution of Equation $(S)$, then $F=f^{-1}$ is a solution of Equation $(S_{0})$. \ES also considers the case of $\phi(z)=z+\lambda$, for a fixed $\lambda\in\C$, which leads to the Abel's equation (Equation $(A)$ is the particular case of $\lambda=1$): 
\begin{equation}
\tag{$A_{0}$}
F(R(z))=F(z)+\lambda.
\end{equation}

Let us note that the case $\phi(z)=z^{d}$, where $d\geq 2$ is an integer, gives rise to the following form of the Böttcher's equation: 
\begin{equation}
\tag{$B_{0}$}
F(R(z))=F(z)^{d}.
\end{equation}

Likewise, if $f$ is an invertible solution of Equation $(B)$, then $F=f^{-1}$ is a solution of Equation $(B_{0})$.

Moreover, E. Schröder is interested in the so-called \textit{analytic iteration problem}, which he formulates in the following way in \cite{Schroder_2}. For a given analytic function $\phi(z)$, find a function $\Phi(w,z)$ of two complex arguments, which is continuous (even analytic) in both variables and such that:
\begin{align}
\label{an_it_pb}
\Phi(w,z)&=\Phi(w-1,\phi(z))\\
\Phi(1,z)&=\phi(z).\nonumber
\end{align}

Even if E. Schröder seems not to explicitly connect this problem to the resolution of the Schröder's equation, it is worth pointing this link out here, as made in \cite{Alexander}. If there exists an invertible analytic solution $F$ to Equation $(S_{0})$, then, a solution of \eqref{an_it_pb} with $\phi=R$ is given by:
\begin{equation}
\label{sol_it_pb_S}
\Phi(w,z)=F^{-1}(s^{w}F(z)).
\end{equation}

An other solution can be given based on an invertible analytic solution $F$ of the Abel's equation $(A_{0})$ by:
\begin{equation}
\label{sol_it_pb_A}
\Phi(w,z)=F^{-1}(F(z)+w\lambda).
\end{equation}

Finally, a solution can be given based on an invertible analytic solution $F$ of the Böttcher's equation $(B_0)$ by:
\begin{equation}
\label{sol_it_pb_B}
\Phi(w,z)=F^{-1}\left(F(z)^{d^{w}}\right).
\end{equation}

Thus, solving Equation $(S), (B)$ or $(A)$ implies solving the analytic iteration problem for the rational fraction $R(z)$. In \cite{Koenigs_1}, G. Koenigs proves the existence of an analytic solution $F$ of Equation $(S_{0})$ around an attracting fixed point $\alpha$, and also (applying his previous reasoning to $F^{-1}$ instead of $F$) around a repelling fixed point $\alpha$. Note that in each case, the solution is unique up to a constant multiplier. Let us consider the attracting case. Recall that $s=R'(\alpha)$. The proof of the author consists in showing that the following function is analytic in a neighbourhood $D$ of $\alpha$.
\begin{equation*}
F(z)=\text{lim}_{n\to +\infty}\frac{R^{n}(z)-\alpha}{s^{n}}.
\end{equation*}

Indeed, if we write
\begin{equation}
\label{def_Fn}
F_{n}(z)=\frac{R^{n}(z)-\alpha}{s^{n}},
\end{equation}
we have :
\begin{align*}
F_{n}(R(z))&=\frac{R^{n+1}(z)-\alpha}{s^{n}}\\
&=sF_{n+1}(z).
\end{align*}
Then, taking the limit when $n$ tends to infinity, we obtain that $F$ is solution of Equation $(S_{0})$.

Now, to prove that $F$ is analytic in a neighbourhood of $\alpha$, G. Koenigs reduces the problem to the convergence of a series of functions $\sum f_{i}(z)$. He then uses a result from G. Darboux which states that if the series $\sum f_{i}(z)$ and $\sum f'_{i}(z)$ both converge uniformly in a disc, then $\sum f_{i}(z)$ converges to an analytic function on this disk. Let us reproduce here a shorter proof we can find in \cite[page 49]{Audin} for example, replacing G. Darboux result by a theorem which states that a locally uniform limit of a sequence of analytic functions in an open set $D$ is analytic on $D$. First, there exist $\sigma\in\R$, such that $s<\sigma<1$, a real number $r>0$ little enough, and $A\in\R^{*}_{+}$, such that $R$ is analytic in the disc $D(\alpha,r)$ centred at $\alpha$ and of radius $r$, and such that for all $z\in D(\alpha,r)$:
\begin{align}
\label{ineq_1}
|R(z)-\alpha|&<\sigma|z-\alpha|,\\
\label{ineq_2}
|R(z)-\alpha-s(z-\alpha)|&<A|z-\alpha|^{2}.
\end{align}

This arises from the Taylor formula applied to $R(z)$ near $z=\alpha$ and the fact that $\alpha$ is an attracting fixed point of $R$. We deduce from the first inequality above that $R(D(\alpha,r))\subset D(\alpha,r)$. Besides, even if this means reducing $r$ we assume that $\alpha$ is the only solution to $R(z)=\alpha$ in $D(\alpha,r)$.

Moreover, we see that for all $z\in D(\alpha,r)$:
\begin{equation}
\label{product_koenigs}
F_{n}(z)=(z-\alpha)\prod_{k=0}^{n-1}\frac{R^{k+1}(z)-\alpha}{s(R^{k}(z)-\alpha)}.
\end{equation}

The fact that $R(D(\alpha,r))\subset D(\alpha,r)$ implies by induction that  $R^{n}(D(\alpha,r))\subset D(\alpha,r)$, for every integer $n\geq 1$. Then, for all $z\in D(\alpha,r)$ we can apply \eqref{ineq_1} to $R^{k-1}(z)$ and \eqref{ineq_2} to $R^{k}(z)$. This gives, for all $z\in D(\alpha,r)\setminus\{\alpha\}$:
\begin{equation*}
\left|\frac{R^{k+1}(z)-\alpha}{s(R^{k}(z)-\alpha)}-1\right|<\frac{A}{|s|}|R^{k}(z)-\alpha|=:u_{k}(z).
\end{equation*}

But it appears that for all $z\in D(\alpha,r)\setminus\{\alpha\}$:
\begin{equation*}
\frac{u_{k+1}(z)}{u_{k}(z)}<\sigma<1.
\end{equation*}

We conclude that the product of \eqref{product_koenigs} converges uniformly in $D(\alpha,r)\setminus\{\alpha\}$. Moreover, for every $n\in\N$, $F_{n}(z)$ admits the finite limit $0$ when $z$ tends to $\alpha$. By (see \cite[7.11 Theorem]{Rudin}), the sequence $(F_{n}(z))$ of analytic functions over $D(\alpha,r)$ converges uniformly on $D(\alpha,r)$. We deduce that $F$ is analytic in $D(\alpha,r)$. Notice that for every integer $n\geq 1$, $\left(R^{n}\right)'(\alpha)=s^{n}$. Furthermore, the uniform convergence of $(F_{n}(z))_{n}$ near $\alpha$ makes the operations of derivation and limit commute in \eqref{def_Fn}. We deduce that $F'(\alpha)=1$. Thus $F$ is locally invertible. This result of G. Koenigs strengthen the connection between the Schröder's equation and the theory of iteration.
\medskip

For his part, P. Fatou points out the connections \eqref{sol_it_pb_S} and \eqref{sol_it_pb_A} between Schröder's and Abel's equations and the analytic iteration problem. Indeed, in \cite{Fatou_1}, P. Fatou remarks that the solution of the Schröder's equation given by G. Koenigs solves this problem. Then, P. Fatou studies the Abel's equation $(A_{0})$ in the case of a rationally indifferent fixed point $\alpha$. Inspired by a previous result from L. Leau \cite{Leau_1,Leau_2}, P. Fatou proves in \cite{Fatou_1} the existence of particular domains (that is open connected sets) $L_{1},\ldots, L_{k}$ such that, for all $j\in\{1,\ldots,k\}$, $\alpha$ belongs to the boundary of $L_{j}$, $R(L_{j})\subset L_{j}$ and the restriction to $L_{j}$ of $R^{n}$ converges uniformly to $\alpha$ on $L_{j}$, when $n$ tends to infinity. Such domains are called the \textbf{petals} of $R$ and the result is called the \textit{Leau-Fatou Flower Theorem}. Note that a similar result is proved by G. Julia in \cite{Julia}. In each petal $L_{j}$, P. Fatou proves the existence of an analytic solution of the Abel's equation $(A)$, and mentions again that this provides a solution to the analytic iteration problem. Let us precise here that Equation $(A)$ is introduced by N. H. Abel in \cite{Abel}. The author remarks that it can be turned into a difference equation. Indeed, if $f$ is an invertible solution of Equation $(A)$, then $F=f^{-1}$ is solution of:
\begin{equation}
\tag{$\tilde{A}$}
F(z+1)=R(F(z)).
\end{equation} 
\medskip

Now, let us consider the Böttcher's equation $(B_{0})$, in the case of a super-attracting fixed point $\alpha$. According to \Ritt \cite{Ritt_Bottcher}, the existence of an analytic solution in the neighbourhood of $\alpha$ is due to L. Böttcher in \cite{Bottcher_1,Bottcher_2,Bottcher_3}. \Ritt provides a short proof of this result in \cite{Ritt_Bottcher}. Let us sketch it here. Even if it means considering, instead of $R(z)$, the conjugation of $R(z)$ with an appropriate Möbius transformation \eqref{def_Mobius}, we may assume that $\alpha=0$. Now, let
\begin{equation}
R(z)=\sum_{n=d}^{+\infty}a_{n}z^{n},
\end{equation} 
where $d\geq 2$. Considering $a_{d}^{1/(d-1)}R(z/a_{d}^{1/(d-1)})$, we may assume that $a_{d}=1$. We know that there exists a disc $D$ around zero in which the only zero of $R(z)$ in $D$ is the origin. Even if this means reducing the radius of $D$, we assume that $R(D)\subset D$. This can be deduced from the Taylor formula and the fact that $0$ is a super-attracting fixed point of $R(z)$. By induction, we find that $R^{n}(D)\subset D$, for every integer $n\geq 1$, and that the origin is the only zero of $R^{n}(z)$ in $D$. Moreover, the origin is a zero of $R(z)$ of order $d$. Hence, by induction, the origin is a zero of $R^{n}(z)$ of order $d^{n}$, for every integer $n\geq 1$. Then, for every integer $n\geq 1$: 
\begin{equation}
\label{frac_Bottcher}
\frac{R^{n}(z)}{z^{d^{n}}}\neq 0, \forall n\in\N, \forall z\in D.
\end{equation}

Hence, for every integer $n\geq 1$, there exists a $d^{n}$-th root of \eqref{frac_Bottcher} on $D$, that is an analytic function $g_{n}(z)$ over $D$ such that
\begin{equation*}
g_{n}(z)^{d^{n}}=\frac{R^{n}(z)}{z^{d^{n}}}, \forall z\in D.
\end{equation*}

Now if we let $h_{n}(z)=zg_{n}(z)$, which is analytic over $D$, we get:
\begin{equation*}
h_{n}(z)^{d^{n}}=R^{n}(z), \forall n\geq 1, \forall z\in D.
\end{equation*}
We can then write: $h_{n}(z)=[R^{n}(z)]^{1/{d^{n}}}$. 

Now, we remark that:
\begin{equation}
\label{prod_hH}
h_{n+1}(z)=z\prod_{i=0}^{n}\left[g_{1}\left(R^{i}(z)\right)\right]^{1/d^{i}},
\end{equation}

Some calculation (see the details in \cite[page 188]{Fatou_1}) prove that the product \eqref{prod_hH} converges uniformly on $D$. We deduce that $h_{n}(z)$ converges uniformly on $D$ to an analytic function $F(z)$ over $D$. Finally, we have:
\begin{align*}
\left[R^{n}(R(z))\right]^{1/d^{n}}&=\left[R^{n+1}(z)\right]^{1/d^{n}}\\
&=\left[\left(R^{n+1}(z)\right)^{1/d^{n+1}}\right]^{d}.
\end{align*}
If we take the limit when $n$ tends to infinity, we obtain that $F$ is solution of $(B_{0})$, which concludes the proof.  
\medskip

Let us return to the interest of \ES for the analytic iteration problem. As said before, the author seems not to relate this question to solutions of the Schröder's or Abel's equation in the manner of \eqref{sol_it_pb_S} and \eqref{sol_it_pb_A}. But the author has another connection in mind. He links the existence of a solution to the analytic iteration problem \eqref{an_it_pb} to the one of a continuous curve which contains the iterates of $\phi(z)$. This consideration of invariant structures of the complex plane with respect to the iteration is at the heart of the theory developed by P. Fatou and G. Julia. However, the study of \ES and G. Koenigs are limited to a neighbourhood of a fixed point. One of the main innovations of the work of P. Fatou and G. Julia is their idea and tools to investigate the whole complex plane, and in fact the \textbf{compactification} $\hat{C}=\C\cup\{\infty\}$ of $\C$, dividing it into zones depending on the behaviour of the sequence of iterates of a rational fraction. Let us note that $\hat{C}$ makes the study of rational fractions a central topic, as they are the only analytic functions over $\hat{C}$. 
In order to present the main discoveries of P. Fatou and G. Julia, let us first introduce some notations and definitions. As before, and for now on, we let $R(z)$ denote a rational fraction of degree at least two, and we let $R^{n}(z)$ denote the $n$-th iterate of $R(z)$. 


Quite natural questions arise when considering a point $z_{0}$ close to a fixed point $\alpha$ and the sequence of iterates $z_{n}=R^{n}(z_{0})$. From a local point of  view, we can wonder if there exists a neighbourhood of $\alpha$ in which $(z_{n})$ converges. From a global point of view, we can study the behaviour of this sequence outside such a neighbourhood, and on its boundary. The work of P. Fatou and G. Julia is about the second point. We may see this as the study of the impact of the point $z_{0}$ and its neighbourhood over the convergence of the sequence $(z_{n})$. To translate this impact, P. Fatou \cite{Fatou_1,Fatou_2} involves the theory of normal families developed by P. Montel \cite{Montel}. We say that a family $\mathcal{F}$ of analytic functions defined on a domain $D$ of $\hat{C}$ is \textbf{normal} over $D$ if from every infinite subsequence of $\mathcal{F}$, we can extract a sub-sequence of $\mathcal{F}$ which converges uniformly locally on $D$ (that is in every compact set of $D$). The link with our subject is given by the application of the \textit{Arzela-Ascoli theorem}  to $\mathcal{F}=\{R^{n}\}_{n}$. Indeed, this states the equivalence for a family $\mathcal{F}$ of continuous functions defined on a domain $D$ of $\hat{C}$ to be normal over $D$ or equicontinuous over $D$. But, by definition, $\{R^{n}\}_{n}$ is \textbf{equicontinuous} over $D$ if for every $z\in D$ and every $\epsilon>0$, there exists $\delta>0$ such that for every $n\in\N$ and every $z_{0}\in D$:
\begin{equation*}
|z-z_{0}|<\delta\Longrightarrow|R^{n}(z)-R^{n}(z_{0})|<\epsilon.
\end{equation*}
Thus, the notion of equicontinuity exactly transcribes the fact that the behaviour of the sequence $(z_{n})$ depends on $z_{0}$. In order to understand the boundary of this property, P. Fatou defines and studies the set of all the points of $\hat{C}$ for which the family $\{R^{n}\}_{n}$ is not normal (that is there exists no neighbourhood of these points in which the family is normal). He denotes $F$ this set, for \textit{Frontière} (french word for \textit{boundary}). This set is nowadays written as $J(R)$ (or $J$), for \textit{the Julia set associated to $R$}, and it is its complement $\hat{C}\setminus J(R)$ that is denoted as $F(R)$ (or $F$), this time because of the initial letter of Fatou, and called the \textit{Fatou set associated with $R$}. At the same time, G. Julia defines the set $E$ consisting of all the repelling fixed points of all the iterates $R^{n}$, $n\in\N$, and studies the derived set $E'$ composed by all the accumulation points of $E$, which coincides with $J$ (see Theorem \ref{thm_FJ_2} below). The study of the set $J$ can provide information on the solutions of functional equations. Let us illustrate this fact with the analytic extension of a solution of the Schröder's equation $(S_{0})$ in a neighbourhood of an attracting fixed point $\alpha$, following a method introduced by P. Fatou \cite{Fatou_3}. Recall that $s=R'(\alpha)$. According to G. Koenigs, there exists a neighbourhood $D_{\alpha}$ of $\alpha$ and an analytic function $F$ on $D_{\alpha}$ such that for all $z\in D_{\alpha}$:
\begin{equation*}
F(R(z))=sF(z).
\end{equation*}

Even if this means reducing $D_{\alpha}$, we assume that $R(D_{\alpha})\subset D_{\alpha}$. This comes from the Taylor formula and the fact that $\alpha$ is attracting.
Let us consider $\tilde{D}=R^{-1}(D_{\alpha})$. The fact that $R(D_{\alpha})\subset D_{\alpha}$ implies that $D_{\alpha}\subset\tilde{D}$. Then, let us define for all $\tilde{z}\in\tilde{D}$:
\begin{equation*}
\tilde{F}(\tilde{z})=\frac{1}{s}F(R(\tilde{z})).
\end{equation*}

We see that 
\begin{equation*}
\tilde{F}(z)=F(z), \forall z\in D_{\alpha}.
\end{equation*}

Moreover, as $R(\tilde{z})\in D_{\alpha}$, for all $z\in\tilde{D}$, we have:
\begin{align*}
\tilde{F}(R(\tilde{z}))&=F(R(\tilde{z})), \text{ } \forall z\in\tilde{D}\\
&=s\tilde{F}(\tilde{z}), \text{ } \forall z\in\tilde{D}.
\end{align*}

Hence, $\tilde{F}$ extends analytically $F$ on $\tilde{D}$ and remains a solution of the Schröder's equation over $\tilde{D}$. Iterating the process, we can extend $F$ to an analytic function $G$ over the domain of attraction $D$ of $\alpha$ (that is the set of all the elements $z$ such that there exists an integer $N$ such that $R^{N}(z)\in D_{\alpha}$), such that $G$ remains a solution of the Schröder's equation over $D$. Thus, the understanding of properties of $R(z)$, namely the nature of its fixed points, can provide information about a solution $F(z)$ of the Schröder's equation $(S_{0})$. 
\medskip

However, Equation $(S)$, the other form of the Schröder's equation, also provides a connection with the theory of P. Fatou and G. Julia, maybe even better in some cases. Indeed, when $\alpha$ is a repelling fixed point of $R$, H. Poincaré \cite{Poincare} proves that the solution $f$ of the Schröder's equation $(S)$ (also called \textit{Poincaré's equation} in this form) can be extended as a meromorphic function over $\C$. Hence, this allows a global study of the structure of the Julia set $J$. As explained in \cite[Theorem 6.3.2]{Beardon}, a computation of the coefficients of the Taylor series expansion of a formal solution $f(z)$ of $(S)$ implies that this series has a positive radius of convergence. Hence, there exists $r>0$ such that $f(z)$ is analytic over the disc $D(\alpha,r)$ centred at $\alpha$ and of radius $r>0$. Then, one can extend $f(z)$ by induction as follows. The analyticity of $f(z)$ over $D(\alpha,r)$ implies that $R(f(z))$ is meromorphic over $D(\alpha,r)$. Hence, by Equation $(S)$, so is $f(sz)$. Therefore, $f(z)$ is meromorphic over $D(\alpha,sr)$. By induction, we find that $f(z)$ is meromorphic over $D(\alpha,s^{n}r)$, for every $n\in\N$. As $|s|>1$, we obtain that $f(z)$ is meromorphic over $\C$.
\medskip

The Schröder's equation is also studied from an algebraic point of view. Indeed, \Ritt establishes in \cite{Ritt} the list of all the differentially algebraic solutions of this equation, when $|s|>1$, that is, in the case of a repelling fixed point of $R(z)$. This is Theorem \ref{thm_BB} of Section \ref{sec_ing_1}.

\subsection{The 80's and the result of \Beck and \Berg}
\label{sec_80}

During the 30's, the interest for functional equations $(S)$, $(B)$ and $(A)$ and the theory of iteration is less vigorous. Nonetheless, let us note the work of C. L. Siegel \cite{Siegel_42} on the existence of a solution to Equation $(S)$ for some indifferent fixed points, and the one of H. Brolin \cite{Brolin} on the structure of the Julia set. The enthusiasm for this subject rises again sixty years later, during the 80's. This is mainly due to the connections the theory of P. Fatou and G. Julia shares with the active area of dynamical systems and fractals, and the possibility to make computational experiences. Indeed, as said before, the iteration of a rational function $R(z)$ gives rise to a dynamical system, which divides $\hat{C}$ into different areas, depending on the concordance or disparity of the local behaviour of the sequence $\{R^{n}(z)\}_{n}$ around a point $z=z_{0}$. The concordance is formalized by the notion of equicontinuous and normal families. The set of points with this local concordance is the Fatou set $F(R)$ and its complement in $\hat{C}$ is the Julia set $J(R)$. Let us assume that $R(z)$ is of degree at least two. Note that $F(R)$ is open in $\hat{C}$ and $J(R)$ is a closed compact subset of $\hat{C}$. Let us mention that $J(R)$ is always non-empty \cite{Fatou_1} and perfect, that is, closed and without any isolated point \cite{Fatou_2, Julia}. Moreover, Julia sets provide lots of examples of fractals. These objects are defined by B. Mandelbrot in the 90's \cite{Mandelbrot}. We refer to \cite{Douady,Chabert,Peitgen} for more details about what follows. The story of fractals actually goes back to the works of B. Riemann and K. Weierstrass and their discoveries of continuous functions with no derivative, at any point. This created a lot of confusion in the mathematical community which had trouble to apprehend such strange objects (see more details in \cite[page 88]{Alexander})! This discomfort was even increased with the work of G. Cantor and his perfect, totally disconnected (that is with all connected components reduced to a point) sets, which questioned the notion of dimension of his time. To deal with the complexity of such objects, the classical topological dimension is replaced by the Hausdorff dimension, which may be a non integral number. Intuitively (see for example \cite{Peitgen}), this notion measures the growth of the number of sets of diameter $\epsilon$ needed to cover the concerned set, when $\epsilon$ tends to zero. The Hausdorff dimension is always greater than or equal to the topological one. B. Mandelbrot defines the fractals as sets for which the Hausdorff dimension is strictly greater than the topological one. For example, the triadic Cantor set has topological dimension 0 and Hausdorff dimension $\log(2)/\log(3)$ : it is a fractal. As said before, lots of Julia sets provide examples of fractals. Julia sets still fuel the current research, with different points of view, namely, investigations on their Hausdorff dimension \cite{Levin,Yang}, their Lebesgue measure \cite{Cheritat}, or their computational complexity \cite{Dudko}. Apart from the structure of the Julia set $J(R)$, there is the question of its variation when the coefficients of $R(z)$ depend on a parameter. This question is raised by P. Fatou in \cite{Fatou_2}. In the case of polynomials of degree two, of the form $R_{c}(z)=z^{2}+c$, where $c\in\C$, there exists a classification of the Julia sets $J(R_{c})$. Note that each polynomial of degree two is conjugated to a polynomial of such a form. This classification is encoded by the Mandelbrot set, defined as the set of the complex numbers $c\in\C$ for which $J(R_{c})$ is connected. Indeed, the dynamic of $R_{c}(z)$ changes as $c$ moves from a cardioid to a disc of the Mandelbrot set. For example (see \cite[paragraph 1.6]{Beardon}), let us focus on the cardioid $C$, and the disc $D$, where $C=u(D(0,1/2)), u(z)=z-z^{2}$, and $D=D(-1,1/4)$. The cardioid $C$ and the disc $D$ are part of the Mandelbrot set (see the figure below). When $c\in C$, the rational fraction $R_{c}$ admits a unique attracting fixed point $\alpha$ and a pair $(u,v)$ such that $R_{c}(u)=v, R_{c}(v)=u$, and $u,v$ are repelling fixed points of $R_{c}^{2}$. Then, when $c$ enters inside $D$, the point $\alpha$ becomes a repelling fixed point of $R_{c}$ and the points $u,v$ become attracting fixed point of $R_{c}^{2}$. As an illustration, the figure below represents the Mandelbrot set and the form of the Julia sets $J(R_{c})$ attached to $R_{c}$ for some of the points $c$ inside and outside the Mandelbrot set. Note that the Mandelbrot set appears to be connected itself \cite{Douady_connexe}.

\begin{figure}[h!]
	\centering
	\includegraphics[width=\textwidth]{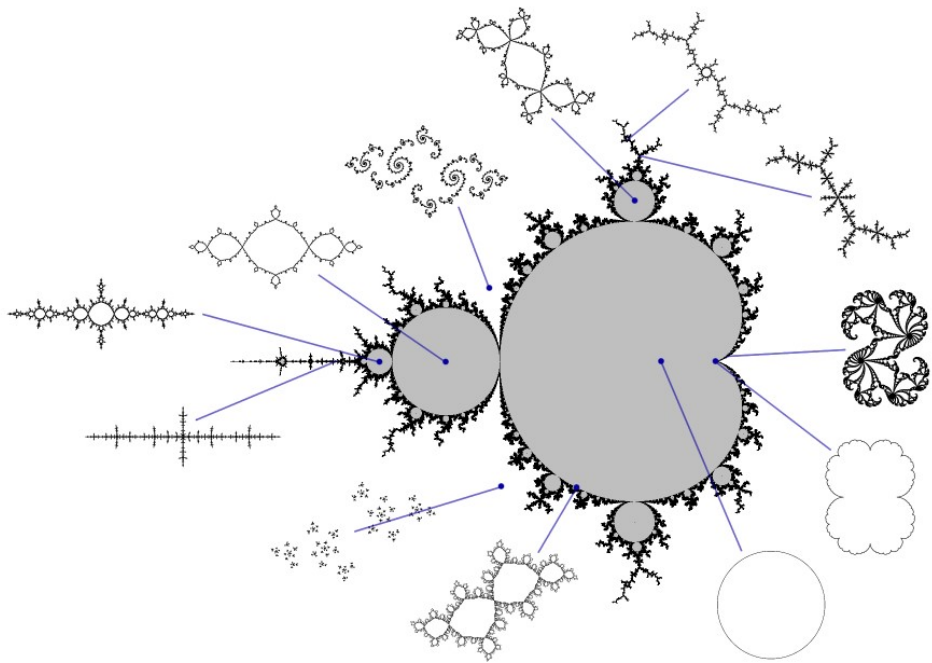}
	\caption{The Mandelbrot set (in grey) and Julia sets $J(R_{c})$ (at the end of the blue lines) for some points $c\in\C$ inside and outside the Mandelbrot set. This image is taken from \cite{Cheritat_image}.}
\end{figure}

Besides, the Mandelbrot set $M$ gives rise to the following interesting transcendental result \cite{N}. A certain conformal map $\Phi(z)$, constructed by A. Douady and J. H. Hubbard \cite{Douady_connexe}, defined on the complement of $M$, admits transcendental values $\Phi(\alpha)$ at each algebraic number $\alpha$ of the complement of $M$. A link with the Böttcher's equation is that $\Phi(c)=f_{c}(c)$, for every $c$ in the complement of $M$, where $f_{c}$ satisfies Equation $(B_{0})$ for $d=2$ and $R(z)=R_{c}(z)$. In the same area, the algebraic aspects of the solutions of Equations $(S), (B)$ and $(A)$ are also investigated by many authors, along with other kinds of functional equations. The first main result of this type is due to 0. Hölder in 1887. Indeed, the author proves that the Euler's Gamma function defined for every $z\in\C$ such that $\text{Re}(z)>0$, by:
\begin{equation*}
\Gamma(z)=\int_{0}^{+\infty} t^{z}e^{-t}\frac{dt}{t},
\end{equation*}
is hypertranscendental over $\C(z)$. As said in \cite{Rubel_survey}, the proof of O. Hölder is based on the following functional equation:
\begin{equation*}
\Gamma(z+1)=z\Gamma(z).
\end{equation*}

This is a non-autonomous version of the Abel's difference equation $(\tilde{A})$, of the form: 
\begin{equation*}
G(z+1)=R(z,G(z)),
\end{equation*}
where $R(z,X)$ is this time a complex rational fraction of two variables. In \cite[Problem 69]{Rubel_2}, L. A. Rubel proposes the study of such a generalised functional equation for the Schröder's equation. This gives rise to further studies. For example, K. Ishizaki \cite{Ishizaki} considers the case where $R(z,X)=a(z)X+b(z)$, with $a(z),b(z)$ rational fractions. The author proves that every transcendental meromorphic solution of the associated generalised Schöder's equation:
\[G(sz)=a(z)G(z)+b(z),\]
with $|s|\notin\{0,1\}$, is hypertranscendental.
\medskip

Concerning the generalised equation of $(B)$, that is:
\begin{equation}
\label{Mahler_eq}
G(z^{d})=R(z,G(z)),
\end{equation}
this is called a $d$-Mahler equation and was introduced by K. Mahler in 1929 in \cite{M1,M2,M3}. A solution $G$ of \eqref{Mahler_eq} is called a $d$-Mahler function. Ke. Nishioka \cite{Keiji_Nishioka} proves that a Mahler function is transcendental if and only if it is not rational (see also \cite{N}). The same author provides in \cite{Ke_Nish} a sufficient condition for the hypertranscendence of Mahler functions of order one. Before this time, K. Mahler proved that the $d$-Mahler function $\sum_{n=0}^{+\infty} z^{d^{n}}$ is hypertranscendental \cite{Mahler}. This result is generalised in \cite{Loxton} with the study by J. H. Loxton and A. J. Van der Poorten of inhomogeneous linear Mahler systems of order one. Such functional results are motivated by a theorem of K. Mahler \cite{M1,M2,M3} which establishes, under some assumptions, the equivalence between the transcendence of a Mahler function $f(z)$ and the one of its value $f(\alpha)$ at a non-zero algebraic number $\alpha$. Thus, results on functional algebraic independence turns into results on algebraic independence of values. More results are obtained using an adapted Galois theory. For short, this theory (see \cite{vdP-Singer} for linear difference equations) associates to a system of linear functional equations an algebraic group, called the \textit{Galois group}, which encodes the algebraic relations between the solutions of the system. If the Galois group is \textit{big enough}, the solutions are algebraically independent. This approach allows K. Nguyen \cite{Nguyen_Mahler} to recover the hypertranscendental result of Ke. Nishioka mentioned above. To illustrate further use of Galois theory for linear Mahler functions, let us mention the result in \cite{Hardouin} which provides sufficient conditions for a Mahler function to be hypertranscendental, and the work of J. Roques \cite{JR2}. Note that these kinds of functional results are still open questions for linear Mahler equations over a function field of positive characteristic (see for example \cite{GF_1,GF_2,GF_these}). Similarly, the Siegel-Shidlovskii theorem \cite{S-S}, which is a kind of analogue of the Mahler's theorem for certain solutions of linear differential equations over $\C(z)$ called $E$-functions, motivates the study of the algebraic independence of solutions of such equations. As for Mahler functions, there is a dichotomy between rational and transcendental $E$-functions. Moreover, the Hrushovski-Feng algorithm \cite{Hrushovski,Feng} computes the Galois group of systems of linear differential equations. More generally, the study of hypertranscendence or algebraic independence of solutions of different types of functional equations is currently a very dynamic area (see for example \cite{Lucia_survol,Hardouin_survol}). This frequently involves the development and use of Galois theories adapted to the different settings. As an illustration, let us mention \cite{Hardouin_Singer} for linear difference equations, and the work of C. Hardouin for $q$-difference systems \cite{Hardouin_hypertransc}, which generalises the work of K. Ishizaki mentioned above, and the development of a general tannakian Galois theory, which in particular apply for $\tau$-difference systems in positive characteristic \cite{Hardouin_Tannakian}, developed by G. Anderson, W. D. Brownawell and M. Papanikolas \cite{ABP1}.
\medskip

Let us go back to Equations $(S), (B)$ and $(A)$. A fruitful link between algebraic properties of solutions of these equations and the iteration theory is given in 1986 by a result of M. Boshernitzan and L. A. Rubel in \cite{Boshernitzan_Rubel}. This states the equivalence for a solution of Equation $(S), (B)$ or $(A)$ to be differentially algebraic and the family $\{R^{n}(z)\}_{n}$ to be \textbf{coherent}. The latter means that all the rational fractions $R^{n}(z)$ satisfy a same algebraic differential equation \eqref{def_DA}. Note that L. A. Rubel asked in \cite[Problem 27]{Rubel_2} if there exists a transcendental entire function whose iterates form a coherent family. This is answered negatively by \Berg in \cite{Bergweiler}.

The theorem of M. Boshernitzan and L. A. Rubel, is one of the ingredients of the proof of the result of P.-G. Becker and W. Bergweiler \cite{Becker_Bergweiler} we are interested in. This statement, reproduced here as Theorem \ref{thm_BB}, explicitly lists all the differentially algebraic solutions of Equations $(S), (B)$ and $(A)$. Let us note that partial results have already been found. We mentioned earlier the result of \Ritt \cite{Ritt}. But we can also indicate the work of F. W. Carroll \cite{Carroll}. The author considers the case where $R(z)$ is a finite Blaschke product with an attracting fixed point (see for example the survey \cite{Blaschke} for the definition and more information about such products). Then, he guarantees that a solution of the Schröder's equation $(S_{0})$ is hypertranscendental. Furthermore, P. Borwein examines the case where $d=2$, and $R(z)=z^{2}+c$, with $c>0$, in the Böttcher's equation $(B)$. Then, the author states that every solution of this equation is hypertranscendental. Moreover, P. Borwein points out that his proof, as well as the one of \Ritt in \cite{Ritt}, shares analogies with the proof of the hypertranscendence of the Gamma function by O. Hölder. Finally, P.-G. Becker and W. Bergweiler themselves previously obtained in \cite{Becker_Bergweiler_0} the list of all the algebraic solutions of the Böttcher's equation $(B)$ when $R(z)$ is conjugated to a polynomial. In fact, this list concerns the solutions of the more general following equation:
\begin{equation}
\label{eq_gen_BB}
f(p(z))=q(f(z)),
\end{equation}
where $p(z),q(z)$ are polynomials of the same degree $d\geq 2$, and the attracting fixed point is $\infty$.

In the same paper \cite{Becker_Bergweiler_0}, the authors conjecture that the transcendental solutions of Equation \eqref{eq_gen_BB} are in fact hypertranscendental. For the Böttcher's equation $(B)$, Theorem \ref{thm_BB}, their further result, implies that this conjecture is satisfied. To conclude, let us mention some recent works. K. D. Nguyen \cite{Nguyen_local_conj} studies systems of $n$ Böttcher's equations $(B)$ for polynomials $R_{1}(z),\ldots, R_{n}(z)$. The author proves a result that links the algebraic independence of some transformations of the solutions $f_{R_{i}}$  associated to each Böttcher's equation to the conjugacy of some iterates of some polynomials among $R_{1}(z),\ldots, R_{n}(z)$. Finally, M. Aschenbrenner and W. Bergweiler prove in \cite{Ashenbrenner_Bergweiler} the hypertranscendence over $\C(z)$ of the iterative logarithm $\text{itlog}(R)$ of a non-linear rational or entire function $R$ with a rationally indifferent fixed point. The function $\text{itlog}(R)$ is the unique formal power series solution $f$ of the equation:
\begin{equation}
\label{Abel_derivee}
f(R(z))=R'(z)f(z),
\end{equation} 
The proof of the authors is similar to the one of Theorem \ref{thm_BB} we present later. Note that the authors also prove that in the case where $R(z)$ is a non-linear entire function, the function $\text{itlog}(R)$ is even hypertranscendental over the ring of entire functions. Note that Equation \eqref{Abel_derivee} is useful to study the iteration of $R(z)$ not only in the petals of $R(z)$, given by the Leau-Fatou Flower theorem mentioned above, but in a neighbourhood of the concerned fixed point.

\section{The complete classification of differentially algebraic solutions of the Schröder's, Böttcher's and Abel's equations}

\subsection{The statement}
\label{sec_Statement}

Before presenting the statement of \Beck and \Berg we are interesting in, let us precisely define its setting. Let $R(z)$ denote a rational fraction with coefficients in $\C$ and of degree at least two. Even if this means replacing $R(z)$ by a conjugate with the appropriate Möbius transformation, we may assume that $0$ is a fixed point of $R(z)$. Let us write $s=R'(0)$. The framework considered is the following (see \cite[Chapitre II]{Fatou_1}):

\begin{enumerate}
	\item 
	
	If $0$ is an attracting, repelling or irrationally indifferent fixed point of $R(z)$, we consider the Schröder's equation $(S)$. The Schröder's equation admits a unique solution of the form
	\begin{equation}
	f(z)=\sum_{n=1}^{+\infty}a_{n}z^{n}, \text{ with } a_{1}=1.
	\end{equation} 
	
If $f$ converges in a neighbourhood of $0$, we say that $f$ is a Schröder function. In the case of an attracting or repelling fixed point, the solution is always convergent, as seen before. But in the case of an irrationally indifferent fixed point, there always exists a formal solution, but not necessarily convergent \cite{Siegel_42}. The question of the convergence of this formal solution is a dynamic area of research \cite{Yoccoz,Geyer}.
	\item 
	If $0$ is a super-attracting fixed point of $R(z)$, there exists an integer $d\geq 2$ such that:
	\begin{equation}
	R(z)=\sum_{n=d}^{+\infty}b_{n}z^{n}, \text{ } b_{d}\neq 0.
	\end{equation}
	Then, we consider the Böttcher's equation $(B)$. For every $a_{1}\in\C$ such that $a_{1}^{1-d}=b_{d}$, the Böttcher's equation admits a unique solution of the form
	\begin{equation}
	f(z)=\sum_{n=1}^{+\infty}a_{n}z^{n}.
	\end{equation} 
	The function $f$ converges in a neighbourhood of $0$, and we say that $f$ is a Böttcher function.
	\item 
	If $0$ is a rationally indifferent fixed point of $R(z)$, even if this means replacing $R(z)$ by an iterate $R^{k}(z)$, we may assume that $s=1$. Indeed, if $s^{k}=1$, we have $(R^{k})'(0)=s^{k}=1$. Let us write:
	\begin{equation}
	R(z)=z+\sum_{n=d}^{+\infty}b_{n}z^{n}, \text{ where } d\geq 2,\text{ } b_{d}\neq 0.
	\end{equation}
	Then, we consider the Abel's equation $(A)$. In each petal given by the Leau-Fatou Flower theorem mentioned earlier, the Abel's equation admits an analytic solution $f(z)$. We say that $f$ is an Abel function.
	
\end{enumerate}

Let us introduce some vocabulary. Let us remind that we say that two analytic functions $S_{1}$ and $S_{2}$ are \textbf{conjugated} to each other via an invertible analytic function $g$ if $S_{1}=g^{-1}S_{2}g$, where functional composition is denoted multiplicatively. If we say that $S_{1}$ and $S_{2}$ are \textbf{conjugated}, with no precision about $g$, we mean that $g$ is a Möbius transformation \eqref{def_Mobius}.
\medskip

We are now able to present the result of \Beck and W. Bergweiler, as it is written in \cite{Becker_Bergweiler}. We give several explanations and comments about Theorem \ref{thm_BB} directly after its statement.

\begin{thm_en}[\Beck and \Berg]
	\label{thm_BB}
Let $f(z)$ be a Schröder, Böttcher or Abel function. Let us assume that $f$ is differentially algebraic. Then, we have the following. 

\begin{enumerate}
	\item 
	If $f$ is a Schröder function, then $0$ is a repelling fixed point of $R(z)$ and $f$ is a Möbius transformation of a function of one of the following forms:
	\begin{enumerate}
		\item
		$\exp(\alpha z^{r})$. In this case, $R(z)$ is conjugated to $z^{d}$ or $z^{-d}$.
		\item 
		$\cos(\alpha z^{r}+\beta)$. In this case, $R(z)$ is conjugated to $T_{d}$ or $-T_{d}$, where $T_{d}$ is the $d$-th Tchebychev polynomial.
		\item 
		$\wp(\alpha z^{r}+\beta)$, $\wp^{2}(\alpha z^{r}+\beta)$, $\wp^{3}(\alpha z^{r}+\beta)$, $\wp'(\alpha z^{r}+\beta)$, where $\wp$ denotes the Weierstrass function,
	\end{enumerate}
where the constant $r$ is a rational number such that the concerned functions are meromorphic over $\C$, $\alpha$ is a non-zero complex number and $\beta$ is a fraction of a period of the concerned function.
\item 
	If $f$ is a Böttcher function, then, $f$ is a Möbius transformation or is a Möbius transformation of a function of one of the following forms:
	\begin{enumerate}
		\item 
		 $\rho z$, where $\rho^{d-1}=1$. In this case, $R(z)$ is conjugated to $z^{d}$.
		 \item 
		 $\rho z+\frac{1}{\rho z}$,  where $\rho^{d-1}=1$. In this case, $R(z)$ is conjugated to $T_{d}$.
		 \item 
		 $\rho z+\frac{1}{\rho z}$, where $\rho^{d-1}=-1$. In this case, $R(z)$ is conjugated to $-T_{d}$.
	\end{enumerate}
\item
The function $f$ is not a Abel function. 
	
\end{enumerate}

In particular, Abel functions are always hypertranscendental.

\end{thm_en}

The meromorphic condition, for the Schröder functions, is due to the fact that, in the case of a repelling fixed point, every solution of Equation $(S)$ extends to a meromorphic function on the complex plane, as seen at the end of Section \ref{sec_birth}. Remark that, even if the function $z\to z^{r}$ is not meromorphic on the complex plane when $r$ is not an integer, this case may happen. There is an example below, with the meromorphic function $\cos(\sqrt{2z})-1$ on the complex plane, where $r=1/2$.
\medskip

We refer to the proof of Theorem \ref{thm_FJ_1}, and references inside, for more details on the rational fraction of the Schröder's equation associated to the Weierstrass function (see also \cite{Ritt_0}).

\medskip
We reproduced in Theorem \ref{thm_BB} the statement as it is written in \cite{Becker_Bergweiler}. The fixed point of the rational fractions $R(z)$ involved in this statement is not necessarily $0$. We clarify this point below.

\begin{enumerate}
	\item 
	If $f$ is a Schröder function
	\begin{enumerate}
		\item
		The rational fractions $R(z)=z^{d}$ or $R(z)=z^{-d}$ are considered at the repelling fixed point $z=1$. To move it at the origin, we have to conjugate $R(z)$ with $L(z)=z-1$. Then we have: $\tilde{R}(z)=(z+1)^{d}-1$, or $\tilde{R}(z)=(z+1)^{-d}-1$, and associated solutions of equation $S_{\tilde{R},0}$ are of the form $\exp(\alpha z^{r})-1$, where $\alpha,r$ are as in the statement of Theorem \ref{thm_BB}. In particular, $\tilde{f}(z)=\exp(z)-1$ is such that $\tilde{f}(0)=0$ and $\tilde{f}'(0)=1$.
		\item 
		The rational fractions $R(z)=T_{d}$ or $R(z)=-T_{d}$ are considered at the repelling fixed point $z=1$. To move it at the origin, we have to conjugate $R(z)$ with $L(z)=z-1$. Then we have: $\tilde{R}(z)=T_{d}(z+1)-1$, or $\tilde{R}(z)=T_{-d}(z+1)-1$, and associated solutions of equation $S_{\tilde{R},0}$ are of the form	$\cos(\alpha z^{r}+\beta)-1$, where $\alpha,\beta,r$ are as in the statement of Theorem \ref{thm_BB}. In particular, $\tilde{f}(z)=\cos(\sqrt{2z})-1$ is such that $\tilde{f}(0)=0$ and $\tilde{f}'(0)=1$.
	\end{enumerate}
	
	\item 
	If $f$ is a Böttcher function,
	\begin{enumerate}
		\item 
		The rational fraction $R(z)=z^{d}$ is considered at the super-attracting fixed point $z=0$.
		\item 
		The rational fraction $R(z)= T_{d}(z)$ is considered at the super-attracting fixed point $z=\infty$. To move it at the origin, we have to conjugate $R(z)$ with $L(z)=1/z$. Then we have: $\tilde{R}(z)=1/(T_{d}(1/z))$, and associated solutions of equation $B_{\tilde{R},0}$ are of the form $1/(\rho z+\frac{1}{\rho z})$, where $\rho^{d-1}=1$.
		\item 
		The rational fraction $R(z)= T_{-d}(z)$ is considered at the super-attracting fixed point $z=\infty$. To move it at the origin, we have to conjugate $R(z)$ with $L(z)=1/z$. Then we have: $\tilde{R}(z)=1/(T_{-d}(1/z))$, and associated solutions of equation $B_{\tilde{R},0}$ are of the form $1/(\rho z+\frac{1}{\rho z})$, where $\rho^{d-1}=-1$.
	\end{enumerate}
		
\end{enumerate}

Let us recall the following definitions. The \textbf{Weierstrass function} $\wp$ is a meromorphic function over $\C$, attached to a \textbf{lattice} $\Lambda$$\subset\C$, that is a discrete subgroup of $\C$ that contains an $\R$-basis of $\C$, and defined by:
\begin{equation*}
\wp(z):=\wp_{\Lambda}(z)=\frac{1}{z^{2}}+\sum_{w\in\Lambda, w\neq 0}\left(\frac{1}{(z-w)^{2}}-\frac{1}{w^{2}}\right).
\end{equation*}

The function $\wp_{\Lambda}$ is \textbf{periodic} with respect to $\Lambda$, that is:
\begin{equation}
\label{Weierstrass_periodique}
\wp_{\Lambda}(z+\omega)=\wp_{\Lambda}(z),\text{ }\forall z\in\mathbb{C}, \forall\omega\in\Lambda.
\end{equation}

Besides, the function $\wp_{\Lambda}$ satisfies the following algebraic differential equation:

\begin{equation}
\wp_{\Lambda}'^{2}=4\wp_{\Lambda}^{3}-g_{2}(\Lambda)\wp_{\Lambda}-g_{3}(\Lambda),
\end{equation}
where $g_{2}(\Lambda), g_{3}(\Lambda)\in\mathbb{C}$ depend on $\Lambda$ and satisfy $g_{2}(\Lambda)^{3}-27g_{3}(\Lambda)^{2}\neq 0$. For more details, see for example \cite{Silverman}.
\medskip

Finally, \textbf{Tchebytchev polynomials} are defined by induction with:
\begin{align*}
T_{0}(X)=1,T_{1}(X)=X,
T_{n+2}(X)=2XT_{n+1}(X)-T_{n}(X),\forall n\geq 2. 
\end{align*}

Let us note that for every integer $n\geq 0$ and every $x\in[-1,1]$: $T_{n}(\cos(x))=\cos(nx)$.

\subsection{First ingredient: Preliminary results around the (hyper)transcendence of solutions of Equations $(S)$ and $(B)$}
\label{sec_ing_1}

\subsubsection{A result of \Ritt}

As mentioned earlier, \Ritt gives in 1926 the list of all the differentially algebraic Schröder functions when $0$ is a repelling fixed point of $R(z)$ \cite{Ritt}. These are exactly the functions listed in the first point of Theorem \ref{thm_BB}. Thus, when the Schröder's equation admits a differentially algebraic solution, then, the considered fixed point of the associated rational fraction is always repelling.

\begin{thm_en}[\Ritt]
	\label{thm_Ritt}
Let us assume that $0$ is a repelling fixed point of a rational fraction $R(z)$ of degree at least two. Let $f$ be a solution of the associated Schröder's equation $(S)$. If $f$ is differentially algebraic, then, $f$ is in the list of the first point of Theorem \ref{thm_BB}.
\end{thm_en}

The proof of this theorem is based on the theories of differentiation and elimination. These techniques allow the author to prove that a solution of Equation $(S)$ is (after an appropriate change of variables) a solution of a Schwarzian differential equation of the following form:
\begin{equation}
\label{eq_Schwarz}
g^{(3)}g'-3/2(g^{(2)})^{2}=A(g)g^{(4)},
\end{equation}
where $A$ is a rational fraction.

The author then uses a previous classification of his own \cite{Ritt_0} for differentially algebraic solutions of \eqref{eq_Schwarz} to deduce that they are of the desired forms. As noted earlier, P. Borwein in \cite{Borwein} points out that the techniques of this proof share analogies with the proof of the hypertranscendence of the Gamma function by O. Hölder.

\subsubsection{A result of \Beck and \Berg}

One year before proving Theorem \ref{thm_BB}, \Beck and \Berg \cite{Becker_Bergweiler_0} provided the list of all the algebraic Böttcher functions, when $R(z)$ is conjugated to a polynomial. It is exactly the functions listed in the second point of Theorem \ref{thm_BB}. Thus, the Böttcher functions are transcendental if and only if they are hypertranscendental. This was a part of a conjecture of the two authors in \cite{Becker_Bergweiler_0}. Let us state this previous result in the case of Böttcher functions. As said before, this applies to the more general equations \eqref{eq_gen_BB}.

\begin{thm_en}[\Beck and \Berg]
	\label{thm_BB_0}
	Let $R(z)$ be a polynomial of degree at least two and let $f$ be a solution of the associated Böttcher's equation $(B)$. Then, $f$ is algebraic if and only if $f$ is in the list of the second point of Theorem \ref{thm_BB}. In particular, all the algebraic Böttcher functions are rational.
\end{thm_en}

The proof of Theorem \ref{thm_BB_0} is based on the analysis of the finite singularities (algebraic branch points) of the Böttcher function. Indeed, the authors prove that if $f$ is an algebraic solution of $(B)$, different from a Möbius transformation, then its local inverse $f^{-1}$ has exactly two such finite singularities. This result, and \cite[Theorem 4.1.2]{Beardon} about exceptional points of a rational fraction (see \cite[Definition 4.1.1]{Beardon}), allow them to find the corresponding forms of $R(z)$. Then, the first remarks of the paper \cite{Becker_Bergweiler_0} provide the Böttcher functions when $R(z)\in\{z^{d},T_{d},-T_{d}\}$. Notice that a Möbius transformation of an algebraic function remains algebraic.

\subsection{Second ingredient: the notion of coherent families}
\label{sec_ing_2}

Remind that a formal power series $f(z)$ is \textbf{differentially algebraic} if there exists a non-zero polynomial $P(z,X_{0},\ldots, X_{n})$ such that $f$ satisfies \eqref{def_DA}. Examples of such functions are given by polynomials, algebraic functions, Bessel functions and classical functions as the exponential, logarithm, cosinus or sinus for example. As said before, the first example of hypertranscendental function does not appear before 1887, with the Euler's Gamma function given by O. Hölder. But \textit{most} entire functions, or analytic functions over a domain of $\C$, are hypertranscendental \cite[Theorem 5]{Rubel_survey}.
\medskip

In \cite[ Problem 26'']{Rubel_2}, L. A. Rubel asks the following question: are there any boundary on the growth of entire differentially algebraic functions ? More precisely, given an entire differentially algebraic function satisfying an $n$-order equation \eqref{def_DA}, do there exist constants $A,\alpha$ such that:
\begin{equation}
\label{q_Rubel}
|f(z)|\leq A \exp^{n}(|z|^{\alpha})?
\end{equation}

L. A. Rubel indicates that a strategy to answer the question negatively should be to construct a function $f$ \textit{big enough} (compared to the exponential) such that all its iterates $f^{k}(z)$ satisfy the same algebraic differential equation. This is the notion of coherent family studied by M. Boshernitzan and L. A. Rubel in \cite{Boshernitzan_Rubel}. Indeed, recall that a family of functions is said to be \textbf{coherent} if all its elements satisfy the same algebraic differential equation \eqref{def_DA}. As an example, let us consider the family $\{cz^{n}, c\in\C, n\in\N\}$. This family is coherent because each of its elements satisfies the following algebraic differential equation:
\begin{equation}
\label{eq_poly_coherent}
zf^{(2)}(z)f(z)+f(z)f'(z)-zf'(z)^{2}=0.
\end{equation}

However, the family of all polynomials with rational coefficients is proved not to be coherent in \cite{Boshernitzan_Rubel}. Before stating the main result of M. Boshernitzan and L. A. Rubel, let us mention two important results concerning coherent families and sketch their proof.

\begin{thm_en}
\label{thm_BR_1}
Let $f$ be an analytic differentially algebraic function. Then, $f$ satisfies an \textit{autonomous} algebraic differential equation. This means that there exist a non-zero polynomial $Q(X_{0},\ldots, X_{n})$, with coefficients in $\C$, independent of $z$, such that:
\begin{equation}
\label{def_DA_auto}
Q(f(z),f'(z),\ldots, f^{(n)}(z))=0.
\end{equation} 
\end{thm_en}  

\begin{proof}
We can find a proof of this well-known result in \cite{Boshernitzan_Rubel}, and a more detailed reasoning in \cite{Rubel_gen}. Let us gather the explanations here. Let $P(z,X_{0},\ldots, X_{n})$ be a non-zero polynomial such that $f(z)$ satisfies Equation \eqref{def_DA}. Without any loss of generality, we can assume that $P$ is irreducible in $\C[z,X_{0},\ldots, X_{n}]$. The goal is then to remove the variable $z$ from the equation. The notion of the resultant of two polynomials will solve the problem. Let us define the following operator on polynomials $S(z,X_{0},\ldots, X_{n})\in\C[z,X_{0},\ldots, X_{n}]$, over $\C[z,X_{0},\ldots, X_{n},X_{n+1}]$:
\begin{equation*}
D: S(z,X_{0},\ldots, X_{n})\longmapsto DS(z,X_{0},\ldots, X_{n+1})=\frac{dS}{dz}(z,X_{0},\ldots, X_{n})+\sum_{k=0}^{n}\frac{dS}{dX_{k}}(z,X_{0},\ldots, X_{n})X_{k+1}.
\end{equation*}

The advantage of this definition is that for all $S\in\C[z,X_{0},\ldots, X_{n}]$ we have:
\begin{equation*}
DS(z,f(z),\ldots, f^{(n+1)}(z))=\frac{d}{dz}[S(z,f(z),\ldots, f^{(n)}(z))]
\end{equation*}

Hence, $DP(z,f(z),\ldots, f^{(n+1)}(z))=0$. Hence, if $R$ is the resultant of the polynomials $P$ and $DP$ with respect to the variable $z$, properties of the resultant guarantee that $R\in\C[X_{0},\ldots, X_{n+1}]$ and that there exist $A,B\in \C[z,X_{0},\ldots, X_{n+1}]$ such that:
\begin{equation}
\label{resultant}
R=AP+B(DP).
\end{equation}
Now, if we specialize \eqref{resultant} at $(z,f(z),\ldots, f^{(n+1)}(z))$, we obtain:
\begin{equation}
\label{auto_eq}
R(f(z),\ldots, f^{(n+1)}(z))=0.
\end{equation}
This provides an autonomous algebraic differential equation for $f$ if we prove that $R$ is a non-zero polynomial. To do so, let us assume by contradiction that $R=0$. Then, $P$ and $DP$ admit a non-constant common factor. As $P$ is irreducible, $P\mid DP$ in $\C[z,X_{0},\ldots, X_{n+1}]$. Let $T\in\C[z,X_{0},\ldots, X_{n+1}]$ such that:
\[(DP)(z,X_{0},\ldots, X_{n+1})=P(z,X_{0},\ldots, X_{n})T(z,X_{0},\ldots, X_{n+1}).\] 
Then, for all $U(z)\in\C[z]$, we have:
\begin{equation}
\label{eq_QU}
(DP)(z,U(z),\ldots, U^{(n+1)}(z))=P(z,U(z),\ldots, U^{(n)}(z))T(z,U(z),\ldots, U^{(n+1)}(z)).
\end{equation}
 Let $U(z)\in\C[z]$. Let us write $Q(z,U(z),\ldots, U^{(n+1)}(z))=\tilde{Q}(z)$, for all polynomial $Q\in\C[z,X_{0},\ldots, X_{n+1}]$. Then, by \eqref{eq_QU} we have:
\begin{equation}
\tilde{P}'(z)=\tilde{(DP)}(z)=\tilde{P}(z)\tilde{T}(z).
\end{equation}
This provides: 
\begin{equation}
\label{constante}
\tilde{P}(z)\in\C. 
\end{equation}
Now, an argument from linear algebra allows us to conclude that $P\in\C$, which is a contradiction. More precisely, if $x_{0},x_{1}\in\C$, there is a surjective morphism of $\C$-vector spaces:
\begin{align}
\phi_{x_{0},x_{1}}:\C_{2n+2}[z]&\longrightarrow \C^{2n+2}\nonumber\\\
U(z)&\longmapsto \left(\begin{pmatrix}
U(x_{0})\\\colon\\U^{(n)}(x_{0})
\end{pmatrix}, \begin{pmatrix}
U(x_{1})\\\colon\\U^{(n)}(x_{1})
\end{pmatrix}\right),
\end{align} 
where $\C_{2n+2}[z]$ denotes the $\C$-vector space of complex polynomials of degree less that or equal to $2n+2$. Then, for all $x_{0}\in\C$ and for all $(z_{0},\ldots,z_{n})\in\C^{n+1}$, there exists $U(z)\in\C_{2n+2}[z]$ such that $(U(x_{0}),\ldots,U^{(n)}(x_{0}))=(z_{0},\ldots,z_{n})$ and $(U(0),\ldots,U^{(n)}(0))=(0,\ldots,0)$. Thus, Equation \eqref{constante} successively applied at $(x_{0},U(x_{0}),\ldots, U^{(n)}(x_{0}))$ and at $(0,U(0),\ldots, U^{(n)}(0))$ implies that
\begin{equation*}
P(x_{0},z_{0},\ldots,z_{n})=P(0,\ldots,0), \text{ } \forall (x_{0},z_{0},\ldots,z_{n})\in\C^{n+2}.
\end{equation*}
This implies that $P=P(0,\ldots,0)\in\C$ and yields a contradiction. Theorem \ref{thm_BR_1} is proved.
\end{proof}

Another useful result about coherent families is that they are stable with respect to many operations. This is the following statement, proved by M. Boshernitzan and L. A. Rubel in \cite{Boshernitzan_Rubel}.

\begin{thm_en}[M. Boshernitzan and L. A. Rubel]
\label{thm_BR_2}
Let $f,g$ be two analytic differentially algebraic functions. Let $P,Q$ be non-zero differential algebraic polynomials providing an autonomous differential algebraic relation for $f$ and $g$ respectively. Let us consider the functions : 
\begin{equation}
\label{stabilite_coherent}
f+g,f-g,f\times g,f/g,fg,fG, fg',
\end{equation}

where $G$ is a primitive of $g$, and the composition of applications is denoted multiplicatively. Then, for every function $h$ in this list, there exists a complex autonomous polynomial $T(X_{0},\ldots,X_{n})$, \textbf{which depends only on} $P$ and $Q$ (\textbf{and not on} $f$ nor $g$) such that:
\begin{equation}
T(h(z),h'(z),\ldots, h^{(n)}(z))=0.
\end{equation}  	

\textbf{In particular}, coherent families are stable under the operations in \eqref{stabilite_coherent}. In other words, if $\mathcal{F}$ and $\mathcal{G}$ are two coherent families, the family $\{f*g\mid f\in\mathcal{F}, g\in\mathcal{G}\}$ is a coherent family, where $*$ is a fixed operation of the set $\{+,-,\times,\div\}$, or the operation of composition. And the two families $\{fG\mid f\in\mathcal{F}, G'\in\mathcal{G}\}$, $\{fg'\mid f\in\mathcal{F}, g\in\mathcal{G}\}$ are coherent families.
\end{thm_en}

\begin{proof}
In order to explain the reasoning, we only sketch the proof for $h=f+g$ (it is similar in the other cases) and the case where $P,Q$ are of order $1$ (to reduce the notations). In other words, we have the equations: $P(f(z),f'(z))=0$ and $Q(g(z),g'(z))=0$, with $P,Q\in\C[X_{0},X_{1}]$ autonomous. As $P$ is autonomous, if we derive the first equation with respect to the variable $z$, we obtain:
\begin{equation}
\label{sep_f}
f''(z)S(P)(f)+\tilde{P}(f(z),f'(z))=0,
\end{equation}
where $\tilde{P}\in\C[X_{0},X_{1}]$, and $S(P)$ is the \textit{separant} of the polynomial $P$. For a polynomial $U(X_{0},\ldots,X_{n})$ in $\C[X_{0},\ldots,X_{n}]$, $S(U)$ is defined by:
\begin{equation*}
S(U)(X_{0},\ldots,X_{n})=\frac{dU}{dX_{n}}(X_{0},\ldots,X_{n}),
\end{equation*}
where the derivation is made with respect to the biggest variable appearing in $U$. For every analytic function $\phi$, we let $S(U)(\phi)=S(U)(\phi(z),\phi'(z),\ldots,\phi^{(n)}(z))$. Similarly, we have:
\begin{equation}
\label{sep_g}
g''(z)S(Q)(g)+\tilde{Q}(g(z),g'(z))=0,
\end{equation}
where $\tilde{Q}\in\C[X_{0},X_{1}]$.

Let us first assume that $S(P)(f)\neq 0$ and $S(Q)(g)\neq 0$. Then, by \eqref{sep_f} and \eqref{sep_g}, there exists a non-zero rational fraction $R\in\C(Z_{0},\ldots,Z_{3})$ such that:
\begin{equation*}
h^{(2)}=R(f,f',g,g'). 
\end{equation*}
Note that $R$ only depends on $P$ and $Q$. But, $g=h-f$ and $g'=h'-f'$. Thus, there exists a non-zero rational fraction $H_{2}\in\C[X_{0},\ldots,X_{3}]$. Such that:
\begin{equation*}
h^{(2)}(z)=H_{2}(f(z),f'(z),h(z),h'(z)).
\end{equation*}
Note that $H_{2}$ only depends on $P$ and $Q$.
We deduce the existence of non-zero rational fractions $H_{3},H_{4}\in\C(X_{0},\ldots,X_{3})$ such that:
\begin{equation*}
\label{eg_H3}
h^{(3)}(z)=H_{3}(f(z),f'(z),h(z),h'(z)).
\end{equation*}

\begin{equation*}
\label{eg_H4}
h^{(4)}(z)=H_{4}(f(z),f'(z),h(z),h'(z)).
\end{equation*}
Note that $H_{3},H_{4}$ only depend on $P$ and $Q$.

Now, let us consider for every $i\in\{2,3,4\}$, $H_{i}(X_{0},\ldots,X_{3})$ as
\begin{equation*}
H_{i}(X_{0},\ldots,X_{3})\in K:=[\C(X_{2},X_{3})](X_{0},X_{1}).
\end{equation*}

The transcendence degree of $K$ over $L:=C(X_{2},X_{3})$ is equal to $2$. Hence, there exists a non-zero monic polynomial $S(X_{2},X_{3})(Y_{0},Y_{1},Y_{2})\in L[Y_{0},Y_{1},Y_{2}]$, where the $Y_{i}$'s are formal variables for $i\in\{0,1,2\}$, such that:
\begin{equation}
\label{annul_Hi}
S(X_{2},X_{3})(H_{2}(X_{0},\ldots,X_{3}),H_{3}(X_{0},\ldots,X_{3}),H_{4}(X_{0},\ldots,X_{3}))=0.
\end{equation} 

Note that the coefficients of $S$ are in $L$. Hence, $S$ only depends on the $H_{i}(X_{0},\ldots,X_{3})$, that is, only on $P$ and $Q$. 

Now, if we let $X_{0}=f(z), X_{1}=f'(z), X_{2}=h(z), X_{3}=h'(z)$ in \eqref{annul_Hi}, we get:
\begin{equation*}
S(h(z),h'(z))(h^{(2)}(z),h^{(3)}(z),h^{(4)}(z))=0.
\end{equation*}

Now, if we formally replace $h^{(j)}(z)$ by a variable $X_{j}$, for every $j\in\{0,\ldots,4\}$, we find a polynomial $T\in\C[X_{0},\ldots, X_{4}]$ such that

\begin{equation*}
T(h(z),h'(z),h^{(2)}(z),h^{(3)}(z),h^{(4)}(z))=0.
\end{equation*}

As, $S\in L[Y_{0},Y_{1},Y_{2}]$ is monic with respect to the variables $Y_{0},Y_{1},Y_{2}$, $T$ is a non-zero polynomial. Finally, as $S$ only depend on $P$ and $Q$, the same is true for $T$. Hence, $T$ provides a non-zero autonomous differential algebraic relation for $h(z)$, which only depends on $P$ and $Q$. 
\medskip

Then, it remains to treat the case where $S(P)(f)=0$, or $S(Q)(g)=0$. 

First, let us notice that there exist integers $k,l\geq 1$ such that the iterates $S^{k}(P)$ and $S^{l}(Q)$ are non-zero constants. Indeed, the separant of a polynomial strictly reduces its total degree. Hence, as polynomials, $S^{k}(P), S^{l}(Q)\neq 0$. Thus, there exist $k_{1}\leq k-1$ and $l_{1}\leq l-1$ such that:
\begin{align}
\label{sep_it_non_nul}
S^{k_{1}}(P)(f)=0 &\text{ and } S^{k_{1}+1}(P)(f)\neq 0\\
S^{l_{1}}(Q)(g)=0 &\text{ and } S^{l_{1}+1}(Q)(g)\neq 0.\nonumber
\end{align}
Note that for all $k_{1}\leq k$ and $l_{1}\leq l$, $S^{k_{1}}(P), S^{l_{1}}(Q)\neq 0$, as polynomials. Besides, we notice that $k,l$ only depend on $P$ and $Q$, but $k_{1}, l_{1}$ depend on $f$ and $g$. The idea is then to apply the first part of the proof to $S^{k_{1}}(P), S^{l_{1}}(Q)$, instead of $P,Q$ respectively, but we need to get rid of the dependence on $f$ and $g$.
\medskip
 
To do so, let us consider all the integers $k_{1}\leq k-1$ and $l_{1}\leq l-1$. As $S^{k_{1}}(P)$ and $S^{l_{1}}(Q)$ are non-zero polynomials, we can consider two formal variables $\tilde{f}_{k_{1}},\tilde{g}_{k_{1}}$ and assume that they formally satisfy \eqref{sep_it_non_nul}, for $k_{1},l_{1}$. We deduce from the first part of the proof that there exists a differential algebraic relation satisfied by the formal variable $h_{k_{1},l_{1}}=\tilde{f}_{k_{1}}+\tilde{g}_{l_{1}}$, which only depends on $S^{k_{1}}(P), S^{l_{1}}(Q)$. Let us note $T_{k_{1},l_{1}}$ the associated non-zero autonomous differential algebraic polynomial, which only depends on $P,Q$ and $k_{1},l_{1}$. Moreover, for every functions $u,v$ which satisfy \eqref{sep_it_non_nul} for some $k_{1},l_{1}$, we have:
\begin{equation*}
T_{k_{1},l_{1}}(u+v)=0.
\end{equation*}

Then, let us note:
\begin{equation*}
T=\prod_{k_{1}\leq k;l_{1}\leq l} T_{k_{1},l_{1}}.
\end{equation*}
Then, $T$ is a non-zero autonomous differential algebraic polynomial which only depends on $P$ and $Q$. Moreover, $T(f+g)=0$ and this concludes the proof.
\end{proof}

Let us notice that if $f$ is an invertible analytic function which is differentially algebraic, then $f^{-1}$ is also differentially algebraic. Indeed, for all $n\in\N$, $\left(f^{-1}\right)^{(n)}(z)$ is a rational fraction of $f(w),f'(w),\ldots, f^{(n)}(w)$, where $w=f^{-1}(z)$. But $f$ is differentially algebraic. Therefore, the family $\{f^{(n)}(w)\}_{n}$ has a finite transcendence degree over $\C(z)$, and so has the family $\left\{\left(f^{-1}\right)^{(n)}(z)\right\}_{n}$.
\medskip

We are finally able to state the main result of M. Boshernitzan and L.A. Rubel \cite{Boshernitzan_Rubel}, which links coherent families to the algebraic properties of Schröder, Böttcher and Abel functions.

\begin{thm_en}[M. Boshernitzan and L.A. Rubel]
	\label{thm_BB_coherent}
	Let $f$ be a Schröder, Böttcher or Abel function. Let $R(z)$ be the associated rational fraction of degree at least two. Then, $f$ is differentially algebraic if and only if the family $\{R^{n}(z)\}_{n\in\N}$ is coherent. 
\end{thm_en}

The proof of Theorem \ref{thm_BB} only uses the direct implication of Theorem \ref{thm_BB_coherent}, that is why we will only reproduce this part of the proof here. 	

\begin{proof}[Proof of the direct implication of Theorem \ref{thm_BB_coherent}]
Let us assume that $f$ is differentially algebraic. The goal is to prove that the family $\{R^{n}(z)\}_{n}$ is coherent. First, let us assume that $f$ is a Schröder function. Then, $g=f^{-1}$ is differentially algebraic and satisfy Equation $(S_{0})$. Thus, we have:
\begin{equation}
R=g^{-1}(s z)g.
\end{equation}
Then, for all $n\in\N$:
\begin{equation*}
R^{n}=g^{-1}(s^{n} z)g.
\end{equation*}
We have seen that the family $\{s^{n} z\}_{n}$ is coherent (see \eqref{eq_poly_coherent}). Hence, by Theorem \ref{thm_BR_2}, $\{R^{n}\}_{n}$ is coherent.
\medskip

If $f$ is a Böttcher function, then, $g=f^{-1}$ is differentially algebraic and satisfy equation $(B_{0})$. Similarly, we obtain
\begin{equation*}
R^{n}=g^{-1}(z^{d^{n}})g.
\end{equation*}
We have seen that the family $\{z^{d^{n}}\}_{n}$ is coherent (see \eqref{eq_poly_coherent}). Hence, by Theorem \ref{thm_BR_2}, $\{R^{n}\}_{n}$ is coherent.
\medskip

Finally, if $f$ is an Abel function, we find:
\begin{equation}
R^{n}=f^{-1}(z+n)f.
\end{equation}
But the family $\{z+n\}_{n}$ is coherent because $\{n\}_{n}$ is (see \eqref{eq_poly_coherent}) and the coherence is stable under addition by Theorem \ref{thm_BR_2}. Hence, again by Theorem \ref{thm_BR_2}, $\{R^{n}\}_{n}$ is coherent.	
\end{proof}

\subsection{Third ingredient: the theory of P. Fatou and G. Julia}
\label{sec_ing_3}

The ingredients from the theory of P. Fatou and G. Julia needed in the proof of Theorem \ref{thm_BB} are stated as Theorems \ref{thm_FJ_1} and \ref{thm_FJ_2} below. Let $R(z)$ be a rational fraction of degree at least two. Let $F(R)$ denote the Fatou set of $R$. Recall that it is the open set of all the elements $z_{0}\in\hat{C}$ for which the family $\{R^{n}(z_{0})\}_{n}$ is normal in a neighbourhood of $z_{0}$. A fundamental result about normal families is the following normality criterion from P. Montel \cite{Montel}.

\begin{thm_en}[P. Montel]
	\label{thm_Montel}
	Let $D$ be a domain in $\hat{C}$. Let $\Omega=\hat{C}\setminus \{0,1,\infty\}$. Then, the family
	\[\mathcal{F}=\{f:D\longrightarrow\Omega\mid f \text{ is analytic over $D$}\}\]
	is normal over $D$.
\end{thm_en}

Among other things, this allows the author to prove various theorems from \'{E}. Picard. One of them states that an analytic function on a punctured neighbourhood of the origin, which admits an essential singularity at $0$, takes all the values of $\hat{C}$, except at most two. 
\medskip

Now, let $J(R)=\hat{C}\setminus F(R)$ denote the Julia set of $R$. Remind that $J(R)$ is a closed compact subset of $\hat{C}$, which is non empty and perfect. The Fatou set, however, can be empty, as we will see in Theorem \ref{thm_FJ_1}. Let us note that the Fatou and Julia sets are compatible with the notions of iteration and conjugation via a Möbius transformation $g$. Indeed, we can prove that:
\begin{equation}
\label{fatou_iteration}
F(R^{n})=F(R), J(R^{n})=J(R), \forall n\in\N.
\end{equation}
Moreover, if $g$ is a Möbius transformation and $S=gRg^{-1}$, we have:
\begin{equation}
\label{fatou_conj}
F(S)=g(F(R)), J(S)=g(J(R)).
\end{equation}

Now, let us discuss the two results used in the proof of Theorem \ref{thm_BB}. The first one deals with the cases where the Fatou set is empty or not.
\begin{thm_en}
	\label{thm_FJ_1}
Let	$R(z)$ be a rational fraction of degree at least two. 
\begin{enumerate}
	\item
Assume that $R(z)$ admits a non repelling fixed point $\alpha$. If $\alpha$ is irrationally indifferent, assume further that the associated Schröder equation $(S)$ admits a convergent solution in a neighbourhood of $\alpha$. Then $F(R)\neq\emptyset$. 
\item
Assume that $R(z)$ is the rational fraction of the Schröder's equation associated with one of the functions 	$\wp(\alpha z^{r}+\beta)$, $\wp^{2}(\alpha z^{r}+\beta)$, $\wp^{3}(\alpha z^{r}+\beta)$, or $\wp'(\alpha z^{r}+\beta)$ which appears in the first point of Theorem \ref{thm_BB}. Then $F(R)=\emptyset$.
\end{enumerate}

\end{thm_en}

The second result used by \Beck and \Berg to prove Theorem \ref{thm_BB} is the following.
\begin{thm_en}
	\label{thm_FJ_2}
Let	$R(z)$ be a rational fraction of degree $d$ at least two. The set of all the repelling fixed points of all the iterates $R^{n}(z)$, $n\in\N$, is dense in $J(R)$.
\end{thm_en}

Let us sketch the proof of Theorem \ref{thm_FJ_1} (see details in \cite{Beardon}). 

\begin{proof}[Proof of Theorem \ref{thm_FJ_1}]
	
	First, we can prove that an attracting or super-attracting fixed point of $R(z)$ belongs to $F(R)$. Indeed, based on the Taylor development of $R$, there exists $\sigma<1$ and a neighbourhood $D$ of the fixed point $\alpha$ such that:
	\[|R(z)-R(\alpha)|\leq \sigma|z-\alpha|, \text{ } \forall z\in D\]
	if $\alpha$ is attracting; and even:
	\[|R(z)-R(\alpha)|\leq \sigma|z-\alpha|^{2}, \text{ }\forall z\in D,\]
	if $\alpha$ is super-attracting.
	Secondly, if $\alpha$ is rationally indifferent, $\alpha\in J(R)$. This is proved by P. Fatou \cite{Fatou_2} and G. Julia \cite{Julia} (see also \cite[Theorem 6.5.1]{Beardon}). But the Leau-Fatou Flower theorem implies the existence of domains $L_{1},\ldots, L_{k}$ called the petals of $R$ such that, for all $j\in\{1,\ldots,k\}$, $\alpha$ belongs to the boundary of $L_{j}$, $R(L_{j})\subset L_{j}$ and the restriction to $L_{j}$ of $R^{n}$ converges uniformly to $\alpha$ on $L_{j}$, when $n$ tends to infinity. This latter property shows that each petal is included in $F(R)$. Hence, $F(R)$ is not empty.
	\medskip
	
	Finally, if $\alpha$ is an irrationally indifferent fixed point of $R(z)$, we can use a theorem stated in \cite{Beardon} which establishes the following equivalence, valid for an indifferent fixed point of $R(z)$:
	\begin{equation}
	\label{equiv_irrat_ind}
    R(z) \text{ is linearizable in a neighbourhood of } \alpha \Leftrightarrow \alpha\in F(R).
	\end{equation}
	We say that $R(z)$ is \textbf{linearizable} in a neighbourhood $D$ of $\alpha$ if there exists an invertible analytic function $g$ over $D$ such that $R$ is locally conjugated via $g$ to a function of the form:
	\[h(z)=\alpha+(z-\alpha)h'(\alpha).\]
	But this is precisely the case for $R$, which is conjugated, via the solution of the associated Schröder's equation $(S)$, to $h(z)=s z$. Then $\alpha\in F(R)$. 
	\medskip
	 
	Besides, for the second part of the theorem, we only sketch the proof for the Weierstrass function $\wp(z)$, following \cite[p. 74]{Beardon}. Recall that this function is \textbf{periodic}, that is, satisfies \eqref{Weierstrass_periodique}. Moreover, this function satisfies the following Schröder's equation:
	\begin{equation}
	\label{eq_S_Weierstrass}
	\wp(2z)=R(\wp(z)),
	\end{equation}	
	for a certain rational fraction $R(z)$. This is the \textit{duplication formula} for elliptic curves (see for example \cite[page 54, page 170]{Silverman} and \cite[Chapter 1]{Lang}). Now, let $D$ be a disc in $\C$ and let $U=\wp^{-1}(D)$. Let $\phi(z)=2z$. Then, $\phi^{n}(U)=2^{n}U$. Hence, for $N$ big enough, $2^{N}U$ will contain a period parallelogram of the lattice $\Lambda$ associated with $\wp$. This means that the values that $\wp$ will take in $2^{N}U$ are exactly the one it takes on $\C$. But it is known, as a consequence of the open mapping theorem, that $\wp(\C)=\hat{C}$. Hence, using \eqref{eq_S_Weierstrass}, we obtain:
	\[R^{N}(D)=R^{N}(\wp(U))=\wp(2^{N}U)=\hat{C}.\] 
	This implies that $\{R^{n}\}_{n}$ is not equicontinuous over $D$. Indeed, the local behaviour of the iterates does not respect the proximity of antecedent points, because $D$ is sent onto the whole Riemann sphere by $R^{N}$. As $D$ is arbitrary, we conclude that $\{R^{n}\}_{n}$ is not equicontinuous over any open subset of $\C$. This implies that $J(R)=\hat{C}$.
\end{proof}

Note that, as quickly mentioned in Section \ref{sec_80}, the case of an irrationally indifferent fixed point $\alpha$ was proved by C. L. Siegel in 1942 for particular cases of $\alpha$, called diophantine fixed points.
\medskip

Now, let us reproduce the proof of Theorem \ref{thm_FJ_2}, as detailed in \cite[Theorem 6.9.2]{Beardon}. First let us recall the following facts. Let $R(z)$ be a rational fraction of degree $d\geq 2$. Then, $R(z)$ is a $d$-fold map of $\hat{C}$ onto itself. That is, for all $w\in\hat{C}$, the equation $R(z)=w$ admits exactly $d$ solutions in $\hat{C}$, when counting multiplicities. Now, we say that $w$ is a critical value of $R(z)$, if there exists $z_{0}\in\hat{C}$ such that $R(z_{0})=w$ and if there is no neighbourhood of $z_{0}$ in which $R(z)$ is injective. But there are only finitely many critical values of $R(z)$ in $\hat{C}$. Indeed, for an element $x\in\hat{C}$, there exists a neighbourhood of $x$ in which $R(z)$ is injective if $R'(z)$ has neither a zero nor a pole at $x$. When $w$ is not a critical value of $R(z)$, there exists exactly $d$ pairwise distinct elements $z_{i}\in\hat{C}$, $i=1,\ldots, d$ such that $R(z_{i})=w$, for every $i\in\{1,\ldots,d\}$.

\begin{proof}[Proof of Theorem \ref{thm_FJ_2}]
The first part of the proof consists in showing that $J(R)$ is contained in the topological closure of the set of all the fixed points of all the $R^{n}(z)$. Then, by \cite[Theorem 9.6.1]{Beardon}, the set of all the non repelling fixed points of all the $R^{n}$ is finite. We deduce that $J(R)$ is contained in the topological closure of the set $\mathcal{P}$ of all the repelling fixed points of all the $R^{n}(z)$. Finally, as, for every integer $n$, each repelling fixed point of $R^{n}(z)$ is contained in the closed set $J(R^{n})=J(R)$, we have that $J(R)$ is the topological closure of $\mathcal{P}$. This gives the result of Theorem \ref{thm_FJ_2}. 
\medskip

It thus remains to show that $J(R)$ is contained in the topological closure of the set of all the fixed points of all the $R^{n}(z)$. It suffices to consider an open set $\mathcal{N}$ of $\hat{C}$ such that $\mathcal{N}\cap J\neq\emptyset$, and prove that $\mathcal{N}$ contains a fixed point of one of the $R^{n}(z)$. Let $w\in\mathcal{N}\cap J\neq\emptyset$. We may assume that $w$ is not a critical value of $R^{2}$. Indeed, the number of such values is finite and we may thus find a non-critical value of $R^{2}$ in a neighbourhood of $w$ in $\mathcal{N}\cap J$. Then, $R^{-2}(w)$ contains at least four distinct points $w_{j}$, $j=1,\ldots, 4$. Indeed, the fact that $d\geq 2$ implies that the degree of $R^{2}$ is more than or equal to four. At least three of these points, say $w_{1}, w_{2}, w_{3}$ are distinct from $w$. Then, we construct three neighbourhoods $\mathcal{N}_{i}$ of $w_{i}$, $i=1,2,3$, whose topological closures are pairwise disjoint, and a neighbourhood $\mathcal{N}_{0}\subset \mathcal{N}$ of $w$, disjoint from every $\mathcal{N}_{i}$ and such that $R^{2}:\mathcal{N}_{i}\longrightarrow \mathcal{N}_{0}$ is a homeomorphism, with reciprocal $S_{i}$, for every $i\in\{1,\ldots, 3\}$. But $\{R^{n}(z)\}_{n}$ is not a normal family on $\mathcal{N}_{0}$, because $w\in\mathcal{N}_{0}\cap J(R)$. Then \cite[Theorem 3.3.6]{Beardon} (which is a corollary of Theorem \ref{thm_Montel}) gives the existence of $z_{0}\in \mathcal{N}_{0}$, $n\geq 1$ and $i\in\{1,\ldots, 3\}$ such that:
\[R^{n}(z_{0})=S_{i}(z_{0}).\]
We deduce that $R^{2+n}(z_{0})=R^{2}(S_{i}(z_{0}))=z_{0}$. Hence $z_{0}$ is a fixed point of an iterate of $R(z)$, contained in $\mathcal{N}$.
\end{proof}

Note that P. Fatou \cite{Fatou_2} and G. Julia \cite{Julia} proved the finiteness of the sets of attracting and rationally indifferent fixed points of all the $R^{n}(z)$. 

\subsection{The proof of Theorem \ref{thm_BB}}

In this section, we reproduce the proof of Theorem \ref{thm_BB} established in \cite{Becker_Bergweiler} by \Beck and \Berg. We give detailed explanations about the way that the three ingredients exposed earlier merge. This is a beautiful illustration of the power of the interactions between distinct domains of mathematics. As mentioned before, the statement of Theorem \ref{thm_BB} belongs to the theory of hypertranscendence of solutions of functional equations. The proof of Theorem \ref{thm_BB} uses previous results of hypertranscendence and coherent families, along with the theory of iteration of a rational fraction and analytic properties of the sets of Fatou and Julia.
\medskip

\begin{proof}[Proof of Theorem \ref{thm_BB}]
We will use the notations of \eqref{sol_eq_conj} and the framework of Section \ref{sec_Statement}. Let $R(z)$ be a rational fraction of degree at least two, which admits $0$ as a fixed point. According to the nature of this fixed point, we let $f$ be either a Schröder solution of $S_{R,0}$ (recall that we assume that this equation admits a convergent solution), a Böttcher solution of $B_{R,0}$, or a Abel solution of $A_{R,0}$. Let us assume that $f$ is differentially algebraic. Our goal is to prove that $f$ cannot be a Abel function and that $f$ is in the list $1$ or $2$ of the statement of Theorem \ref{thm_BB}, depending on whether $f$ is a Schröder or Böttcher function.
\medskip

Let us immediately remark that Theorem \ref{thm_BB_coherent} guarantees that the family $\{R^{n}(z)\}_{n}$ is coherent. To begin with, if $0$ is repelling (that is, $f$ is a Schröder function), we can apply Theorem \ref{thm_Ritt} and get that $f$ is in the list of the first point of Theorem \ref{thm_BB}. 
\medskip

Now, let us deal with the case where $0$ is not repelling. The goal is to reduce this case to the repelling one. Let us first notice that Theorem \ref{thm_FJ_1} implies that $F(R)\neq\emptyset$. Moreover, by Theorem \ref{thm_FJ_2}, the set of all the repelling fixed points of all the iterates $R^{n}(z)$ are dense in $J(R)$. As $J(R)$ is always a non-empty set \cite{Fatou_1}, we deduce that there exist $w\in\hat{C}$ and $k\in\N$ such that $w$ is a repelling fixed point of $R^{k}(z)$. We may then consider a Schröder solution $\phi$  of $S_{R^{k},w}$. But $\{\left(R^{k}\right)^{n}\}_{n}$ is coherent, as a subfamily of the coherent family $\{R^{n}(z)\}_{n}$. Hence, by Theorem \ref{thm_BB_coherent}, $\phi$ is differentially algebraic. 
\medskip

Then, Theorem \ref{thm_Ritt}, implies that there exists a Möbius transformation $g(z)$ such that $g^{-1}R^{k}g$ is one of the rational fractions which appear in the first point of Theorem \ref{thm_BB} (with $d^{k}$ instead of $d$). Let us note $S(z)=g^{-1}R^{k}g$. By \eqref{fatou_iteration} and \eqref{fatou_conj}, the fact that $F(R)\neq\emptyset$ implies that $F(S)\neq\emptyset$. 
\medskip

Hence, by Theorem \ref{thm_FJ_1}, $S(z)\in\{z^{d^{k}}, z^{-d^{k}}, T_{d^{k}}, -T_{d^{k}}\}$. Besides, as $0$ is a non repelling fixed point of $R(z)$, as noticed in Section \ref{sec_ing_3}, $g^{-1}(0)$ is a non repelling fixed point of $S(z)$. But $z^{-d^{k}}$ does not admit any such fixed points, and all the non repelling fixed points of $z^{d^{k}}$, $T_{d^{k}}$ or $-T_{d^{k}}$ are super-attracting (a computation shows that a fixed point of a Tchebychev polynomial distinct from $\infty$ is repelling and that $\infty$ is a super-attracting fixed point). Hence, $S(z)\in\{z^{d^{k}}, T_{d^{k}}, -T_{d^{k}}\}$, and $g^{-1}(0)$ is a super-attracting fixed point of one of these three polynomials. Note that, by Section \ref{sec_ing_3} again, this implies that $0$ is a super-attracting fixed point of $R^{k}(z)$ and $R(z)$. 
\medskip

We deduce that $f$ is a solution of $B_{R,0}$. Iterating this equation, we see that $f$ is a solution of $B_{R^{k},0}$. Now, by \eqref{sol_eq_conj}, $g^{-1}f$ is a solution $\psi$ of $B_{S,g^{-1}(0)}$. But the fact that $S(z)\in\{z^{d^{k}}, T_{d^{k}}, -T_{d^{k}}\}$ guarantees, via Theorem \ref{thm_BB_0} that $\psi$ is algebraic and is a Möbius transformation or a Möbius transformation of one of the functions of the list of the second point of Theorem \ref{thm_BB}. As $f=g\psi$, we deduce that $f$ is a Möbius transformation or a Möbius transformation (composition with the Möbius transformation $g$) of a function of the list of the second point of Theorem \ref{thm_BB}. 
\medskip

Hence, $f(z)$ is solution of $B_{M,0}$, where $M(z)$ is of the form of the rational fractions appearing in the second point of Theorem \ref{thm_BB}. As $f(z)$ is solution of $B_{R,0}$ and $B_{M,0}$, we have $R=M$ and $R(z)$ is of the form of the rational fractions appearing in the second point of Theorem \ref{thm_BB}. 
\medskip

To conclude, we have proved that $f$ is either a Schröder or Böttcher function (thus $f$ is not a Abel function). Moreover, we have proved that, when $f$ is a Schröder function, $f$ is of the form described in the first point of Theorem \ref{thm_BB}, and  when $f$ is a Böttcher function, $f$ is of the form described in the second point of Theorem \ref{thm_BB}. Theorem \ref{thm_BB} is thus proved.

\end{proof}

\textbf{Acknowledgement.} The author would like to warmly thank Alin Bostan, Lucia Di Vizio, and Kilian Raschel for their enthusiastic proposal to write this paper. She would also like to thank the detailed reports of the reviewers which allowed her to improve the fluidity, clarity and precision of this survey. Finally, she would like to thank Lucia Di Vizio for her support and her reviews and corrections of this paper.

\bibliographystyle{smfplain}
\bibliography{Equations_Schroder_Bottcher_Abel_bis}
	
	\end{document}